\documentclass[11pt]{article}
\usepackage{amsfonts,amsmath,amsxtra}
\usepackage{latexsym}
\usepackage{amssymb}

\def\hybrid{\topmargin 0pt      \oddsidemargin 0pt
        \headheight 0pt \headsep 0pt
        \textwidth 16.5cm
        \textheight 23cm
        \marginparwidth 0.0in
        \parskip 5pt plus 1pt   \jot = 1.5ex}
\catcode`\@=11
\def\marginnote#1{}
\newcount\hour
\newcount\minute
\newtoks\amorpm
\hour=\time\divide\hour by60 \minute=\time{\multiply\hour by60
\global\advance\minute by-\hour}
\edef\standardtime{{\ifnum\hour<12 \global\amorpm={am}%
        \else\global\amorpm={pm}\advance\hour by-12 \fi
        \ifnum\hour=0 \hour=12 \fi
      \number\hour:\ifnum\minute<10 0\fi\number\minute\the\amorpm}}
\edef\militarytime{\number\hour:\ifnum\minute<10 0\fi\number\minute}

\def\draftlabel#1{{\@bsphack\if@filesw {\let\thepage\relax
   \xdef\@gtempa{\write\@auxout{\string
      \newlabel{#1}{{\@currentlabel}{\thepage}}}}}\@gtempa
   \if@nobreak \ifvmode\nobreak\fi\fi\fi\@esphack}
        \gdef\@eqnlabel{#1}}
\def\@eqnlabel{}
\def\@vacuum{}
\def\draftmarginnote#1{\marginpar{\raggedright\scriptsize\tt#1}}

\def\draft{\oddsidemargin -0.1truein
        \def\@oddfoot{\sl preliminary draft \hfil
        \rm\thepage\hfil\sl\today\quad\militarytime}
        \let\@evenfoot\@oddfoot \overfullrule 3pt
        \let\label=\draftlabel
        \let\marginnote=\draftmarginnote
\def\@eqnnum{{\rm (\theequation)}
\rlap{\kern\marginparsep\tt\@eqnlabel}%
\global\let\@eqnlabel\@vacuum}  }
\newfont{\Bbbb}{msbm7 scaled 1\@ptsize00}
\newcommand{\zs}{\raise-1pt\hbox{$\mbox{\Bbbb Z}$}}

scaled 1\@ptsize00
\font\sevenmsa=msam6 %scaled 1\@ptsize00
\newfam\msafam
\textfont\msafam=\sevenmsa
\def\hexnumber@#1{\ifnum#1<10 \number#1\else
\ifnum#1=10 A\else\ifnum#1=11 B\else\ifnum#1=12 C\else \ifnum#1=13
D\else\ifnum#1=14 E\else\ifnum#1=15 F\fi\fi\fi\fi\fi\fi\fi}
\def\msa@{\hexnumber@\msafam}
\def\llcorner{\delimiter"4\msa@78\msa@78 }
\def\lrcorner{\delimiter"5\msa@79\msa@79 }
\mathchardef\blacktriangleright="3\msa@49
\mathchardef\blacktriangleleft="3\msa@4A \font\tenmsb=msbm10 scaled
1\@ptsize00
\newfam\msbfam
\textfont\msbfam=\tenmsb \scriptfont\msbfam=\tenmsb
\newdimen\Squaresize \Squaresize=14pt
\newdimen\Thickness \Thickness=0.5pt

\def\Square#1{\hbox{\vrule width \Thickness
   \vbox to \Squaresize{\hrule height \Thickness\vss
      \hbox to \Squaresize{\hss#1\hss}
   \vss\hrule height\Thickness}
\unskip\vrule width \Thickness} \kern-\Thickness}

\def\Vsquare#1{\vbox{\Square{$#1$}}\kern-\Thickness}

\def\numberbysection{\@addtoreset{equation}{section}
        \def\theequation{\thesection.\arabic{equation}}}
\numberbysection

\renewcommand{\theequation}{\thesection.\arabic{equation}}
\def\titlepage{\@restonecolfalse\if@twocolumn\@restonecoltrue\onecolumn
     \else \newpage \fi \thispagestyle{empty}\c@page\z@
        \def\thefootnote{\fnsymbol{footnote}} }

\def\endtitlepage{\if@restonecol\twocolumn \else  \fi
        \def\thefootnote{\arabic{footnote}}
        \setcounter{footnote}{0}}  %\c@footnote\z@ }
\relax

\hybrid
\makeatletter
\newdimen\normalarrayskip            % skip between lines
\newdimen\minarrayskip               % minimal skip between lines
\normalarrayskip\baselineskip \minarrayskip\jot
\newif\ifold             \oldtrue            \def\new{\oldfalse}
\def\arraymode{\ifold\relax\else\displaystyle\fi}%mode of array enrties
\def\eqnumphantom{\phantom{(\theequation)}} % ight phantom in eqnarray
\def\@arrayskip{\ifold\baselineskip\z@\lineskip\z@
     \else
     \baselineskip\minarrayskip\lineskip1\baselineskip\fi}

\def\@arrayclassz{\ifcase \@lastchclass \@acolampacol \or
\@ampacol \or \or \or \@addamp \or
   \@acolampacol \or \@firstampfalse \@acol \fi
\edef\@preamble{\@preamble
  \ifcase \@chnum
     \hfil$\relax\arraymode\@sharp$\hfil
     \or $\relax\arraymode\@sharp$\hfil
     \or \hfil$\relax\arraymode\@sharp$\fi}}

\def\@array[#1]#2{\setbox\@arstrutbox=\hbox{\vrule
     height\arraystretch \ht\strutbox
     depth\arraystretch \dp\strutbox
width\z@}\@mkpream{#2}\edef\@preamble{\halign \noexpand\@halignto
\bgroup \tabskip\z@ \@arstrut \@preamble \tabskip\z@ \cr}%
\let\@startpbox\@@startpbox \let\@endpbox\@@endpbox
  \if #1t\vtop \else \if#1b\vbox \else \vcenter \fi\fi
  \bgroup \let\par\relax
  \let\@sharp##\let\protect\relax
  \@arrayskip\@preamble}
\def\eqnarray{\stepcounter{equation}%
              \let\@currentlabel=\theequation
              \global\@eqnswtrue
              \global\@eqcnt\z@
              \tabskip\@centering              %formulae  centering
              \let\\=\@eqncr
              $$%
            \halign to \displaywidth  \bgroup
             \eqnumphantom \@eqnsel
      \hskip\@centering                               %right tab%
    $\displaystyle  \tabskip\z@ {##}$%
    &\global\@eqcnt\@ne \hskip 2\arraycolsep
         $ \displaystyle  \arraymode{##}$\hfil
    &\global\@eqcnt\tw@ \hskip 2\arraycolsep
         $\displaystyle\tabskip\z@{##}$\hfil
         \tabskip\@centering
    &{##}\tabskip\z@\cr}
\makeatother
\newcommand{\RR}{{\mathbb{R}}}

\newcommand{\R}{\mathbb{R}}
\newcommand{\C}{\mathbb{C}}

\def\IC{\mathbb{C}}

\def\IF{\mathbb{F}}

\def\IP{\mathbb{P}}
\def\IQ{\mathbb{Q}}
\def\IR{\mathbb{R}}
\def\IZ{\mathbb{Z}}
\def\CD {\mathcal{D}}
\def\CE {\mathcal{E}}
\def\CF {\mathcal{F}}
\def\CG {\mathcal{G}}
\def\CH {\mathcal{H}}

\def\CL {\mathcal{L}}
\def\CM {\mathcal{M}}
\def\CN {\mathcal{N}}
\def\CO {\mathcal{O}}

\def\CV {\mathcal{V}}

\def\a {{\alpha}}

\def\pr {\partial}
\def\apr {\overline {\partial }}
\def\jb{\bar{j}}
\def\kb{\bar{k}}
\def\lb{\bar{l}}

\def\wb {\bar{w}}
\def\zb {\bar{z}}

\def\Tr{{\rm Tr}}

\newtheorem{te}{Theorem}[section]%Usage:\begin{te}Statement\end{te}

\newtheorem{prop}{Proposition}[section]
\newtheorem{lem}{Lemma}[section]

\newtheorem{rem}{Remark}[section]

\newcommand\bqa{\begin{eqnarray}}
\newcommand\eqa{\end{eqnarray}}
\def\be{\begin{eqnarray}\new\begin{array}{cc}}
\def\ee{\end{array}\end{eqnarray}}

\def\beq{\begin{equation}}
\def\eeq{\end{equation}}
\def\bse{\begin{subequations}}                %%%SUBEQUATIONS
\def\ese{\end{subequations}}
\def\bp{\begin{pmatrix}}
\def\ep{\end{pmatrix}}

\def\stack#1#2{\raise0.7pt\hbox{$\mathrel{\mathop{#2}\limits^{#1}}$}}
\def\tr{\triangleright}
\def\tl{\triangleleft}
\def\sem{\mathsurround=0pt \raise1pt
\hbox{$\scriptscriptstyle>\!\!$}\:\!\!\tl}
\def\mes{\mathsurround=0pt \tr\!\:\!\raise0.8pt
\hbox{$\scriptscriptstyle\!\!<$}\,}
\def\]{\mathsurround=0pt ]\raise-2pt\hbox{$_\ast$}}
\def\<{\langle}
\def\>{\rangle}

\def\CO{{\cal O}}

\def\CH{\mathcal{H}}

\def\we{\raise-1pt\hbox{$\,\stackrel{\wedge}{,}\,$}}
\def\tr{{\rm tr}\,}
\def\Tr{{\rm Tr}\,}
\def\pr {\partial}

\newcounter{pac}[section]

\newcounter{pacc}[subsection]

0

\title{\bf Archimedean $L$-factors and \\
Topological Field Theories I}
\begin{document}
\author{Anton Gerasimov, Dimitri Lebedev, and Sergey Oblezin}
\date{}
\maketitle

\renewcommand{\abstractname}{}

\begin{abstract}
\noindent {\bf Abstract}.  We propose
a functional integral representation for Archimedean
$L$-factors given by products  of $\Gamma$-functions.
The corresponding functional integral arises in the
description of type A equivariant topological linear sigma model on a disk.
The functional integral representation 
provides in particular an interpretation of the $\Gamma$-function 
as an equivariant symplectic volume of an  infinite-dimensional space   
of holomorphic maps of the  disk to $\IC$. This 
should be considered as a mirror-dual to the
classical Euler integral representation of the $\Gamma$-function.
We give an analogous functional integral representation
of $q$-deformed $\Gamma$-functions using  a
three-dimensional equivariant topological linear sigma model on a
handlebody. In general, upon proper ultra-violent  
completion, the topological sigma model considered
on a particular class of three-dimensional spaces  with  a 
compact K\"{a}hler target space provides a  quantum field theory
description of  a $K$-theory version of Gromov-Witten
invariants.

\end{abstract}
\vspace{5 mm}

\section*{Introduction}

Archimedean local $L$-factors were introduced to simplify functional
equations of global $L$-functions. From the point of view of arithmetic
geometry these factors complete the Euler product
representation of  global $L$-factors by taking into account Archimedean
places of the compactified spectrum of global fields.
A known construction of non-Archimedean local $L$-factors is rather transparent
and uses  characteristic polynomials of the image of  the Frobenius
homomorphism in  finite-dimensional representations  of
 the local Weil-Deligne  group closely related to the local
Galois group. On the other hand, Archimedean $L$-factors
are expressed through products of
$\Gamma$-functions and thus are analytic objects avoiding simple
algebraic interpretation. Moreover,  Archimedean  Weil-Deligne groups
are rather mysterious objects in comparison with  their non-Archimedean
counterparts. For instance, in the case of the field of complex
numbers the corresponding Galois group is trivial while  the Weil-Deligne group is
isomorphic to the multiplicative group $\IC^*$ of complex numbers.
In many cases the multiplicative group $\IC^*$  effectively plays the role of
the  Galois group for  complex numbers \cite{De}, \cite {Be}. Thus, for
instance, by  analogy with the action of the  Galois group ${\rm
Gal}(\overline{\IF}_p/\IF_p)$ on  \'{e}tale cohomology of schemes
over $\overline{\IF}_p$,
the multiplicative group $\IC^*$ acts on  the complexified cohomology of
compact non-singular complex algebraic varieties  providing the
standard Hodge decomposition. The Archimedean Weil-Deligne groups 
also play an important
role in a formulation of  Archimedean local Langlands
correspondence (see e.g. \cite{ABV}).

In a series of papers \cite{GLO2}, \cite{GLO3}, \cite{GLO4} we approach
the problem of a proper interpretation of Archimedean $L$-factors
and in particular of Archimedean Weil-Deligne group
using our previous results on relations between
quantum integrable systems
and representation theory as well as results of
\cite{Gi1}, \cite{Gi2}, \cite {Gi3}, \cite {GiL} on quantum
cohomology and quantum $K$-theory.  In \cite{GLO2} we propose
an  explicit description of the Archimedean local Hecke algebras
in terms of intertwining integral operators.
These integral operators are instances of the Baxter
operator playing an important role in
the construction of explicit solutions of quantum
integrable systems. We demonstrate that,
by  analogy with the generators of non-Archimedean Hecke algebras,
 the common eigenfunctions of the Baxter operators 
are given by ``class one'' Whittaker functions. Remarkably  the corresponding
eigenvalues are precisely  Archimedean $L$-factors in complete
analogy with a non-Archimedean case. 
Here the ``class one''  condition means that  the Whittaker function
rapidly decreases outside some domain and is an  Archimedean counterpart of
the non-Archimedean ``class one'' Whittaker function from \cite{Sh}, \cite{CS}.

In \cite{Gi4} Givental introduced new integral representations of
the Whittaker functions arsing. These representation 
arise naturally  in the description of the
generating functions of the equivariant Gromov-Witten (GW) invariants
of the flag spaces.  In \cite{GiL} a $K$-theory analog of the
GW-invariants was considered and it
was shown that the  generating function of $K$-theory  GW-invariants
are given by $q$-deformed  Whittaker functions.
However one should note that  the
Whittaker functions and its $q$-deformations  considered in the framework of
Gromov-Witten theory are not ``class one'' Whittaker
functions and thus do not have direct arithmetic interpretation.
In \cite{GKLO} we show that  Baxter integral operators appear  as
an essential ingredient of the Givental
integral representations of the Whittaker functions.

The relevance of  the Baxter integral operators
both to the description of Archimedean Hecke algebra and Gromov-Witten
invariants implies a relation between these two
subjects. In  \cite{GLO4} we
confirm these expectations. We construct explicit expressions of $q$-deformed
``class one'' Whittaker functions, $q$-analogs of Baxter operators,
 $q$-analogs of Archimedean $L$-factors and  
relate these constructions to counting problems on the moduli space of holomorphic
maps of $\IP^1$ to the flag manifolds. In particular
$q$-deformed $L$-factors and  $q$-deformed
``class one'' Whittaker functions
corresponding to $(\ell+1)$-dimensional representations of the
Archimedean Weil-Deligne group $\IC^*$ are given by characters of
$\IC^*\times GL_{\ell+1}(\IC)$  modules realized as cohomology
of holomorphic line bundles on the moduli space of holomorphic
maps of $\IP^1$. Thus for $q$-deformed $L$-factor the  corresponding
$\IC^*\times GL_{\ell+1}(\IC)$-module is given by
a limit of the space of degree $d$ holomorphic polynomial maps
of $\IC$ to $\IC^{\ell+1}$  when $d\to \infty$. Here  $\IC^*$ acts
on $\IC$ by  multiplication  and $ GL_{\ell+1}(\IC)$ acts  on $\IC^{\ell+1}$
through the standard representation.
We  stress in \cite{GLO4} that the constructed $q$-deformed local
$L$-factors/Whittaker functions   interpolate
Archimedean and non-Archimedean
local $L$-factors/Whittaker functions.
Non-Archimedean $L$-factor
associated to an  $(\ell+1)$-dimensional representation $V=\IC^{\ell+1}$ 
of the Weil-Deligne group can be identified with a  character of a direct
sum of symmetric  powers of $V$.
Similarly, according to \cite{Sh}, \cite{CS},
$p$-adic ``class one'' Whittaker functions
for a reductive algebraic group $G$ are  given by  characters of  irreducible
finite-dimensional representations of the Langlands dual group
${}^LG_0$ (Shintani-Casselman-Shalika formula).
Shintani-Casselman-Shalika (SCS) formula essentially uses
homomorphims of the local Weil-Deligne group (simple modification of the local
Galois group ${\rm Gal}(\overline{\IQ}_p/\IQ_p))$
to  ${}^LG_0$ and thus provides an arithmetic construction of $p$-adic  Whittaker
function.   Remarkably,  the representation \cite{GLO4}
of  the  $q$-deformed Whittaker function  as
a character of a  $\IC^*\times GL_{\ell+1}(\IC)$-module  reduces to
SCS-formula in the appropriate limit. The the same limit 
$q$-deformed  $L$-factor reduces to  non-Archimedean $L$-factors.
Although the interpolation and  limiting procedure allow to relate
 constructions of local Archimedean and non-Archimedean $L$-factors
a direct analytic construction of the Archimedean $L$-factors given by the product of
$\Gamma$-functions was missing.

In this note we provide  an explicit  construction
of  Archimedean $L$-factors as functional integrals
in equivariant  type A topological linear sigma model on a disk.
We remark  in  \cite{GLO4} that  the construction of
the ``class one''  $q$-deformed Whittaker function implies
 consideration of holomorphic maps of  a disk to flag spaces
 in contrast with \cite{Gi4},  \cite{GiL} where holomorphic
maps  of $\IP^1$ into flag spaces were considered. The type A
topological sigma model considered below is precisely
of the kind that describes GW-invariants.
Specifically,  we consider $S^1\times U_{\ell+1}$-equivariant
type A topological linear sigma model on $D=\{z\,|\,|z|\leq 1\}$
with the target space $V=\IC^{\ell+1}$ deformed by a boundary
observable.  The target space here is a
finite-dimensional space $V$ appearing
in the standard construction of an Archimedean $L$-factor in terms of
a finite-dimensional representation of the Archimedean Weil-Deligne group $\IC^*$.
The group  $S^1$ acts by rotations on  $D$ and $U_{\ell+1}$ acts
in $\IC^{\ell+1}$ via the standard representation.
The functional integral representation for local $L$-factors
is given in Theorem \ref{mainTh}. 
Let us stress that the underlying action of
the considered topological theory is quadratic and
the functional integral allows mathematical formulation
using $\zeta$-function regularization \cite{RS}, \cite{H}.

The functional integral representation of local
Archimedean $L$-factors proposed  in this note
can be interpreted as a calculation of
an equivariant symplectic volume of the space of holomorphic maps  of
a disk $D$ into  complex vector spaces. In particular classical
$\Gamma$-function coincides with an equivariant volume of the space of
holomorphic maps of $D$ to $\IC$.  It appears that the infinite
dimensional symplectic geometry is a proper framework for a mysterious
geometry over  Archimedean places.    

The interpretation in terms of equivariant symplectic volumes 
provides a natural way to deform $\Gamma$-function. 
According to standard arguments
in quantum mechanics equivariant symplectic volume
can be interpreted as a  classical limit of the
partition function of the quantum system obtained by a quantization of the
symplectic manifold. We apply this reasoning to the space of
holomorphic maps of $D$ to $\IC$  and obtain a canonical quantization
of $\Gamma$-functions and therefore of local Archimedean $L$-factors.
Fortunately, thus constructed quantum analogs  
of $\Gamma$-function and local Archimedean $L$-factors
coincide with  $q$-deformations of their classical counterparts
\cite{GLO3}, \cite{GLO4}. For these $q$-deformed/quantum local Archimedean $L$-factors
we provide a representation in terms of the functional integral
in  three-dimensional equivariant topological linear sigma model on a
handlebody $S^1\times D$.  We also argue that upon a proper
ultra-violent completion the constructed
three-dimensional topological sigma model
on a special type of compact three-dimensional spaces with
a K\"{a}hler target space  provides a description of quantum $K$-theory invariants
\cite{Gi3}, \cite{GiL}.

Let us finally refer to the  various previous discussions  of  hidden structures behind
local  Archime-dean $L$-factors related with our proposal.
 The point of view that is most close to the approach
of this note was advocated by Deninger
\cite{D1}, \cite{D2} (see also \cite{M}). Namely,
Deninger  proposed an interpretation of  local Archimedean $L$-factors
as  regularized infinite-dimensional  determinants.
Corresponding infinite-dimensional vector spaces
can be identified with  cohomology groups of algebraic manifold with coefficients
in an Archimedean analog of the rings
introduced by Fontaine \cite{F} for a description of the $p$-adic
Hodge periods. Let us note also
that  equivariant symplectic volumes of the space of maps of a disk into
symplectic manifolds were  discussed in \cite{Gi1} in connection with
Gromov-Witten theory.  Also Vojta \cite{Vo} proposed a detailed set of analogies
 between number theory  and  value distribution
 theory of holomorphic functions  due to Nevanlinna.
Finally long ago Deligne  considered an analogy
between schemes over $\IQ_p$ allowing good integer models and
families of complex manifolds over a disk  \cite{De}.
This analogy motivated the development of the theory of mixed Hodge
structures.

The plan of the paper is as follows. In Section 1 we recall two  basic
constructions of local $L$-factors. In Section 2 we
prove the main statement of the paper (Theorem \ref{mainTh})  by identifying
 Archimedean $L$-factors with a particular functional integral
in the equivariant type A topological linear sigma model on a disk.
In Section 3 we provide interpretation of the functional integral
as an equivariant symplectic volume of the
space of holomorphic maps of the disk into complex vector space.
In Section 4 we consider  a three-dimensional
topological sigma model calculating quantum $K$ theory invariants
introduced in \cite{Gi3}, \cite{GiL}. We demonstrate that the functional
integral in the three-dimensional equivariant topological
linear sigma model on the handle-body reproduces $q$-version
of Archimedean $L$-factor introduced in \cite{GLO3}, \cite{GLO4}.
In Section 5 we conclude with some general remarks
and discuss further directions of research.
In Appendix  the standard facts about $\zeta$-regularization
of infinite-dimensional Gaussian integrals are collected.

{\em Acknowledgments}: The research was supported by
Grants RFBR-08-01-00931-a, 09-01-93108-NCNIL-a.
AG was  also partly supported by Science Foundation Ireland grant.

\section{Two constructions of   $L$-factors}

To put our results in the right perspective
we briefly recall in this Section
some standard facts
on  local and global $L$-functions (see e.g.
\cite{Se},  \cite{T}, \cite{ABV}).
This Section does not
contain new results and can be skipped by a connoisseur.
For standard definitions  and elementary
detailed discussion  see for example \cite{Bu}, \cite{L}.

Two essentially different constructions of $L$-functions are known.
The first one uses automorphic representations 
of  reductive algebraic groups. In this construction one expresses $L$-functions
in terms of a spectrum of the corresponding Hecke algebra.
In  the case  of a local field $F$ we consider representations of the group $G(F)$
 and  in the case of a global field $F$ the  group  $G$ is defined
over adeles $\mathbb{A}_F$. The second
construction of $L$-factors is based on arithmetic properties of the
base field $F$ and uses  homomorphims of its Galois group (more
precisely  Weil-Deligne group) into the Langlands dual group ${}^LG$.
Here the Langlands dual group ${}^LG$ for a field $F$
is defined as an extension of the Galois group ${\rm
  Gal}(\overline{F}/F)$ by a dual algebraic group
${}^LG_0(\overline{F})$. The 
Langlands correspondence claims that  both automorphic and arithmetic
constructions lead to the same set of $L$-functions.

Let us describe these two constructions of $L$-functions for
the base field $F=\IQ$  in more details. We start with the automorphic construction of
$L$-functions.   Let $\mathbb{A}$ be the adele ring of  $\mathbb{Q}$ and $G$ be a
reductive algebraic group. An automorphic representation $\pi$ of
$G(\mathbb{A})$ is an irreducible representation entering the decomposition
of the left representation of $G(\mathbb{A})$  in $L^2(G(\mathbb{A})/G(\mathbb{K}))$,
where $G(\mathbb{K})$ is a  maximal compact subgroup in $G(\mathbb{A})$.
The representation $\pi$ can be characterized by
an automorphic form $\phi$ such that it is an eigenfunction of  any element of the
global Hecke algebra $\CH(G(\mathbb{A}))$. The global Hecke algebra
has natural structure of a  product $\CH(G(\mathbb{A}))=(\otimes_p
\CH_p)\otimes \CH_{\infty}$ of local non-Archimedean Hecke
algebras $\CH_p=\CH(G(\mathbb{Q}_p), G(\mathbb{Z}_p))$ for each
prime $p$ and an Archimedean Hecke algebra
$\CH_{\infty}=\CH(G(\mathbb{R}),K)$, where $K$ is a maximal compact
subgroup in $G(\RR)$. The local Hecke algebra $\CH_p$
is isomorphic to the algebra of characters 
of finite-dimensional representations of  a Langlands dual group
${}^LG$, (e.g. for ${}^LG_0$ being
 $A_{\ell}$, $B_{\ell}$,  $C_{\ell}$,  $D_{\ell}$ the  duals are
 $A_{\ell}$, $C_{\ell}$,  $B_{\ell}$,  $D_{\ell}$ respectively).
For each unramified representation  of $G(\mathbb{Q}_p)$ one can
define an action of $\CH_p$ such that the  automorphic form
$\phi$  is a common eigenfunction of all elements of  $\CH_p$
for all  primes $p$ and thus defines  a set of homomorphisms
$\CH_p\to \mathbb{C}$. Identifying local Hecke algebras with the
algebra of characters  of finite-dimensional representations
 of ${}^LG$ one can describe this set of
homomorphisms as a  set  of conjugacy classes $g_p$  in ${}^LG$.
Given a complex finite-dimensional representation
$\rho_V: {}^LG\to GL(V,\IC)$ one can
construct  $L$-function corresponding to the automorphic form
$\phi$ as the  Euler product
\bqa\label{autoL}
L(s,\phi,\rho_V)=\mathop{{\prod}'}_{p}
L_p(s,\phi,\rho_V)=\mathop{{\prod}'}_{p} \,
\det_{V} (1-\rho_V(g_p)\, p^{-s})^{-1},
\eqa
where $\prod_p'$ is a product over the primes $p$ such that
the corresponding representation of $G(\mathbb{Q}_p)$ is not ramified.
It is natural to complete the product by including  local
$L$-factors corresponding to Archimedean and ramified places.
We omit a discussion of $L$-factors for  ramified representations
and consider the construction of Archimedean $L$-factors.
For  Archimedean places the
Hecke eigenfunction property is  usually replaced  by the eigenfunction property
with respect to a ring of invariant
differential operators on $G(\mathbb{R})$. The corresponding
eigenvalues are described by a conjugacy class  $t_{\infty}$ in the Lie
 algebra ${}^L\mathfrak{g}_0={\rm Lie}({}^LG_0)$. The Archimedean
 $L$-factor  is then given by a product of $\Gamma$-functions
 \be\label{archL}
L_{\IR}(s,\phi,\rho_V)=\prod_{j=1}^{\ell+1}
\,\,\pi^{-\frac{s-\alpha_j}{2}}\,
\Gamma\Big(\frac{s-\alpha_j}{2}\Big)\,
=\,\det_V\,\,\,\pi^{-\frac{s-\rho_V(t_{\infty})}{2}}\,
\Gamma\Big(\frac{s-\rho_V(t_{\infty})}{2}\Big),
\ee
where $\rho:\,{}^LG_0\to GL_{\ell+1}(\IC)$ and 
$\rho_V(t_{\infty})={\rm diag}(\alpha_1,\ldots
\alpha_{\ell+1})$.  In \cite{GLO2} we demonstrate
that  the eigenvalue property
with respect to  the ring of invariant
differential operators on $G(\mathbb{R})$
can be equivalently replaced by the eigenvalue property with respect
to a set of integral operators. These operators 
generate Archimedean Hecke algebra
$\CH(G(\IR),K)$ where $K$ is a  maximal  compact subgroup of $G(\IR)$. 
Elements $\CH(G(\IR),K)$ are  $K$-biinvariant functions on 
$G(\IR)$ and the structure of the algebra is given by a convolution. 
This completely restores the 
symmetry between Archimedean and non-Archimedean places. The local
Archimedean $L$-factors \eqref{archL}  correspond to the real place of
$\IQ$. For more general global
fields complex Archimedean places appear and the corresponding
complex local $L$-factors are given by  analogous expressions
 \bqa\label{archLCone}
L_{\IC}(s,\phi,\rho_V)=\,\det_V\,\,\,(2\pi)^{-(s-\rho_V(t_{\infty}))}\,
\Gamma\,\,(s-\rho_V(t_{\infty})).
\eqa
  The global $L$-function can be defined as an analytic
continuation of the product
\be\label{completeL}
\Lambda(s,\phi,\rho)=L(s,\phi,\rho)L_{\infty}(s,\phi,\rho),
\ee
where $L_{\infty}(s)$ is a product of local factors corresponding to
Archimedean places. Global $L$-functions
should satisfy a  functional equation
\be\label{FuncEq}
\Lambda(1-s,\phi,\rho)=\epsilon(s,\phi,\rho)
\Lambda(s,\phi_{\pi^{\vee}},\rho^{\vee}),
\ee
where the $\epsilon$-factor is of the form
$\epsilon(s,\phi,\rho)=A\,B^s$,\,$A\in\IC^*$, $B\in \IR_+$ 
 and $\pi^{\vee}$, $\rho^{\vee}$
are dual to $\pi$, $\rho$.

Now let us describe another construction of the local $L$-functions
based on  arithmetic properties of a base field $F$.
To discuss both Archimedean and non-Archimedean cases  in the same
setting  one should  introduce a notion of the Weil-Deligne group
generalizing  the standard notion of the Galois group. 
Let $\bar{F}$ be an algebraic closure of a local number field $F$.
The  Weil group $W_{F}$ of $F$
should satisfy  the following properties.
First, there should exist a homomorphism with a dense image in the natural
topology on ${\rm Gal}(\bar{F}/F)$
$$
\phi:\,\, W_F\rightarrow {\rm Gal}(\bar{F}/F).
$$
Second, for any Galois extension $E$ of $F$, there should be
an inclusion in the abelianzation $W^{ab}_E=W_E/[W_E,\,W_E]$
of $W_E$
$$
r:  E^*\rightarrow W^{ab}_E,
$$
such that the composition
$$
E^*\rightarrow {\rm Gal}^{ab}(\bar{F}/E),
$$
is a basic homomorphism  of abelian class field theory.

In the case of the local field $F=\IQ_p$ this boils down to the
following  construction. Let $\mathbb{F}_p$ be a residue field
of $\IQ_p$. Then we have the extension 
\be\nonumber
0\longrightarrow I_p\longrightarrow {\rm Gal}(\bar{\IQ}_p/\IQ_p)
\longrightarrow {\rm Gal}(\bar{\IF}_p/\IF_p)\longrightarrow 0,
\ee
where $I_p$ is the inertia group and 
${\rm Gal}(\bar{\IF}_p/\IF_p)$ is a pro-finite completion of the
cyclic groups with the  generator given by Frobenius homomorphism $Fr_p:\,x\to
x^p$. The Weyl group $W_{\IQ_p}$
 is then a subgroup of ${\rm   Gal}(\bar{\IQ}_p/\IQ_p)$ consisting
of the elements such that their  image in ${\rm
  Gal}(\bar{\IF}_p/\IF_p)$ is an integer  power of Frobenius
$Fr_p$. Maximal abelian subgroup of $W_{\IQ_p}$ is canonically
identified with $\IQ^*_p$.

For the field $F=\IC$ of complex numbers,
 the Galois group  ${\rm Gal}(\overline{\IC}/\IC)$ is trivial
and $W_{\IC}=\IC^*$. The homomorphism $\phi$ is trivial and $r$ is the identity map.
For real numbers $F=\IR$, the Galois group
${\rm Gal}(\IC/\IR)=\IZ_2$ is generated by complex conjugation
$Fr_{\infty}$ and the Weil  group $W_{\IR}=\IC^*\cup j\IC^*$
is generated  by a copy of  $\C^*$ and an element $j$, subjected  to the relations:
\be
j x j^{-1}=\overline{x}, \qquad  j^2=-1\in \IC^*,
\ee
with the maps
\be
\phi:\,\,W_{\IR}\rightarrow {\rm Gal}(\IC/\IR),\qquad
\phi(x)=1,\,\,\,\,\,\phi(jx)=Fr_{\infty},\quad x\in \IC,
\ee
\be
r:\,\,\IR^* \rightarrow W^{ab}_{\R},\qquad r(x)=x.
\ee
Note that $W_{\IR}$ is non-abelian and
its abelianization is $W^{ab}_{\IR}=W_{\IR}/[W_{\IR},W_{\IR}]=\IR^*$
where we identify ${\rm Gal}(\IC/\IR)=\{\pm 1\}$. Note that for
Archimedean places thus defined Weil group coincides with the
Weil-Deligne group.

Let ${}^LG$ be a Langlands  group dual to the reductive algebraic
group $G$ used in the automorphic  construction
of $L$-factors above. 
Now with any continuous  homomorphism of the Weil group in
${}^LG$   and a finite-dimensional representation of ${}^LG$
such that the image of the Weil group is semisimple 
one can associate an $L$-factor. Let us fix a
 complex finite-dimensional  representation $\rho_V: {}^LG\to
GL(V,\IC)$ and consider  its composition with a  homomorphism
$\pi_p:W_{\IQ_p}\to {}^LG$ such that the image of the inertia group is
trivial. Then the corresponding
local non-Archimedean $L$-factor is given by
\be\label{autoLtwo}
L_p(s,\pi,\rho_V)=\det_{V} (1-\rho_V(g_p)\, p^{-s})^{-1},
\ee
and coincides with the one defined using the local  Hecke algebra
\eqref{autoL} under appropriate identification of the parameters.

For  Archimedean fields we consider the special case when the
image of the Weil group under a homomorphism $W\to {}^LG$
is abelian and the compact subgroup is in the
kernel. Thus in the case of $F=\IC$ we have the multiplicative group $\IR_{>0}$
and in the case of $F=\IR$ we have $\IR^*=\IR_{>0}\times \{\pm 1\}$.
The corresponding local Archimedean $L$-factors are given by
\be\label{archLC}
L_{\IC}(s,\phi,\rho_V)=\,\det_V\,\,\,(2\pi)^{-(s -\Lambda)}\,
\Gamma\,\,(s-\Lambda),
\ee
\be\label{archLL}
L_{\IR}(s,\phi,\rho_V,Fr_{\infty}=+1)
=\,\det_V\,\,\,\pi^{-\frac{s-\Lambda}{2}}\,
\Gamma\Big(\frac{s-\Lambda}{2}\Big),
\ee
\be \label{archLLL}
L_{\IR}(s,\phi,\rho_V,Fr_{\infty}=-1)
=\,\det_V\,\,\,\pi^{-\frac{s-\Lambda}{2}}\,
\Gamma\Big(\frac{s+1-\Lambda}{2}\Big),
\ee
where $\Lambda$ is an image of the generator of $\IR_{>0}$ and
$Fr_{\infty}=\pm 1$ describes the  action of the Frobenius
homomorphism $Fr_{\infty}\in {\rm Gal}(\IC/\IR)$ in $V$.

\begin{rem} \label{GenLfac} Local Archimedean $L$ factors are
  introduced in such a way that the
  completed global $L$ function \eqref{completeL} satisfies  the
  functional equation \eqref{FuncEq}. This leaves a freedom to multiply
the local  Archimedean $L$-factor by a function  of the form $A\,
B^s$, $A\in \IC^*$, $B\in \IR_{> 0}$ and  is compensated
by a freedom to redefine  $\epsilon$-factor in \eqref{FuncEq}.
\end{rem}

\section{$\Gamma$-function via 2d Topological Field Theory}

In this section we provide a functional integral representation
of a product of $\Gamma$-functions using
 the Feynman path integral formulation of
a two-dimensional  topological field theory.
This leads to a functional integral representation of 
Archimedean $L$-factors. The two-dimensional theory involved is
$S^1\times U_{\ell+1}$-equivariant version of the
type $A$ topological sigma model
on the disk $D$ with the target space $V=\IC^{\ell+1}$.
Here $U_{\ell+1}$ acts on $\IC^{\ell+1}$ via standard
representation   and $S^1$ acts by rotations $\alpha:z\to
ze^{\imath \alpha}$ on the disk $D=\{z\in \IC|\,|z|\leq 1\}$.

\subsection{Type $A$ Topological sigma-models}

We start by recalling the standard construction of a topological sigma
model  associated with a K\"{a}hler manifold
with non-negative first Chern class.
For general discussion of the two-dimensional topological
sigma models see \cite{W1}, \cite{W2}, \cite{AM} \cite{CMR} and reference therein.

Let $X$ be a compact K\"{a}hler manifold of a complex dimension $(\ell+1)$
and $\CM(\Sigma,X)=Map(\Sigma,X)$ be the space of maps $\Phi:\,\Sigma\to X$ of a
compact Riemann surface $\Sigma$ to $X$.
Let $(z,\zb)$ be local complex coordinates on  $\Sigma$.
We pick a Hermitian metric $h$ on $\Sigma$
and denote $\sqrt{h}\,d^2z$ the corresponding measure on $\Sigma$.
The complex structure defines a decomposition $d=\pr+\apr$,
$\pr=dz\,\pr_z$, $\apr=d\zb\,\pr_{\zb}$
of the differential $d$ acting on the differential forms on $\Sigma$.
Let $K$ and $\bar{K}$ be canonical and anti-canonical
bundles over $\Sigma$.
Let $\omega$ and $g$ be a K\"{a}hler form   and the K\"{a}hler
metric on $X$ and  $T_{\IC}X=T^{1,0}X\oplus T^{0,1}X$ be a
decomposition of the complexified tangent bundle of
$X$. We denote local complex coordinates
 on $X$ by $(\varphi^j$,$\varphi^{\jb})$ . Locally
Levi-Civita  connection $\Gamma$ and
the corresponding Riemann  curvature tensor  $R$
are given by
\be
\Gamma_{jk}^i=g^{i\bar{n}}\pr_{j}g_{k\bar{n}},\qquad
R_{i\bar{j}k\bar{l}}=g_{m\bar{j}}\pr_{\bar{l}}\Gamma^m_{ik}.
\ee

Now let us specify the field content of a type $A$ topological
sigma model. We define commuting fields $F$ and $\bar{F}$ as
sections of $K\otimes \Phi^*(T^{0,1})$ and of $\bar{K}\otimes
\Phi^*(T^{1,0})$ correspondingly. The anticommuting  fields  $\chi$,
$\bar{\chi}$ are sections  of the bundles $\Phi^*(\Pi T^{1,0}X)$,
$\Phi^*(\Pi T^{0,1})$ and anticommuting fields  $\psi$, $\bar{\psi}$ are sections
of the bundles $K\otimes \Phi^*(\Pi T^{0,1})$,  $\bar{K}\otimes
\Phi^*(\Pi T^{1,0})$. Here  $\Pi \CE$ denotes the vector bundle $\CE$ with
the reverse parity of the fibres. Metrics $g$ on $X$ and $h$ on
$\Sigma$ induce a Hermitian paring $\<\,,\,\>$ on
the space of sections of the considered
bundles. Thus for example in local coordinates we have
\be\label{metrex}
\<\chi,\chi\>=\sum_{j=1}^{\ell+1}
g_{i\bar{j}}\,\bar{\chi}^{\jb}\,\chi^{i},\qquad
\<F,F\>=\sum_{j=1}^{\ell+1}
h^{z\zb}g_{i\bar{j}}\bar{F}_z^{\jb}\,F_{\zb}^{i}.
\ee
The action functional of the type A topological sigma model is given by
\be\label{Gsigma}
S_{\Sigma}(t,t^*)=S_0(t)+S_{top}(t^*),
\ee
\be\nonumber
S_0(t)=\int_{\Sigma}d^2z\sqrt{h}\,\Big(\imath \<\CF,\apr
\varphi\>+\imath\<\bar{\CF},\pr \bar{\varphi}\>
+\imath \<\bar{\psi}, D\bar{\chi}\>+\imath\<\psi,\bar{D}\chi\>+
t\<\CF\,\CF\>
+t\<\bar{\psi}, R(\psi,\bar{\chi})\chi\>\Big),
\ee
\be\label{topact}
S_{top}=t^*\int_\Sigma\,\Phi^*(\omega),
\ee
where $\Phi^*(\omega)$ is a pull back of the  K\"{a}hler form $\omega$ on
$X$ along the map $\Phi: \Sigma \to X$. The fields $\CF$, $\bar{\CF}$  and the
derivatives $D_z$, $\bar{D}_{\zb}$ are  locally given by
\be\label{defone}
\CF^j_{\zb}=F^j_{\zb}+\Gamma^j_{kl}\chi^k\psi^l_{\zb},\qquad
\CF^{\jb}_z=F^{\jb}_z+\Gamma^{\jb}_{\kb \lb}\chi^{\kb}\psi^{\lb}_z,
\ee
\be\label{deftwo}
(D_z\bar{\chi})^{\jb}=
\pr_z\bar{\chi}^{\jb}+\Gamma_{\bar{k} \bar{l}}^{\jb}(\varphi)
\pr_z\bar{\varphi}^{\bar{k}}
\bar{\chi}^{\bar{l}},\qquad
(\bar{D}_{\zb}\chi)^{j}=\pr_{\zb}\chi^{j}+\Gamma_{k l}^{j}(\varphi)\pr_{\zb}\varphi^{k}
\chi^{l}.
\ee
Here and in the following we imply the  summation over repeating
indexes.

The action functional \eqref{Gsigma} is  invariant with respect to the
transformations $A\to A+\epsilon \delta A$
where $\epsilon$ is an anticommuting  parameter
and the action of $\delta$ on the fields is
\be\label{transformA}
\delta \varphi=\chi,\qquad \delta \chi=0,\qquad
\delta \bar{\psi}= \bar{F},\qquad \delta \bar{F}=0,
\ee
$$
\delta \bar{\varphi}=  \bar{\chi},\qquad \delta \bar{\chi}=0,
\qquad \delta \psi= F,\qquad
 \delta F=0,
$$
and  we have
\be\label{nilpotence}
 \delta^2=0.
\ee
The action $S_0(t)$ in \eqref{Gsigma}  can be written in the
following form
\be\label{QcomA}
S_0(t)=\int_{\Sigma}\,\,d^2z\sqrt{h}\,\delta\CV(t),
\ee
where
$$\CV(t)=\<\psi,\frac{1}{2}t\CF+\imath \apr\varphi\>+
\<\bar{\psi},\frac{1}{2}t\bar{\CF}+\imath \pr\bar{\varphi}\>.
$$
The action \eqref{Gsigma} is obviously $\delta$-invariant. Indeed,
$S_0(t)$ is $\delta$-invariant due to
\eqref{QcomA} and \eqref{nilpotence} and
 the second term $S_{top}$ is a topological invariant equal to the
degree of the map $\Phi: \Sigma\to X$ and thus invariant under
arbitrary infinitesimal deformations of the fields.

The transformation \eqref{transformA} can be interpreted as an
action of the de Rham differential in the infinite-dimensional setting.
Consider an odd vector bundle $\CE\to \CM(\Sigma,X)$  over the space
$\CM(\Sigma,X)$ of  maps
  $\Phi:\Sigma \to X$ with a fiber isomorphic to the direct sum
  $\left(K\otimes \Phi^*(T^{0,1}X)\right)\oplus \\ \left(
\overline{K}\otimes \Phi^*(T^{1,0}X)\right)$.
 Then \eqref{transformA}  defines the action of $\delta$ 
which coincides with action of 
the  de Rham   differential on differential forms on the infinite-dimensional
supermanifold $\CE$ where we imply $\chi:=\delta\varphi$ and $F:=\delta\psi$.

In the topological fields theories we are interested in
calculating  the Feynman path integrals with the action \eqref{Gsigma}
of $\delta$-closed functionals $A$ of the fields
$\varphi,\bar{\varphi},F,\bar{F},\psi,\bar{\psi}$, $\chi$, $\bar{\chi}$
\be\label{funcint}
\Big<A\Big>_{\Sigma}=\int [\CD\varphi]\,
[\CD\bar{\varphi}]\,[\CD
F]\,[\CD\bar{F}]\,[\CD\psi]\,[\CD\bar{\psi}]\,
[\CD\chi]\,[\CD\bar{\chi}]\,
\,e^{-S_{\Sigma}(t,t^*)}\,A(\varphi,\bar{\varphi},F,\bar{F},\psi,\bar{\psi},
\chi,\bar{\chi}).
\ee
In general a functional integral is mathematically  not well-defined
and is usually handled using the standard renormalization routine of  Quantum
Field Theory. In the following we consider
particular correlation functions in a special kind of type $A$ topological
sigma-model when all  calculations are mathematically well-defined.
Here we proceed formally implying that the integral
\eqref{funcint} is well-defined.

We would like to consider only correlation functions
\eqref{funcint} of (the product of)  observables $A$ which are $\delta$-closed.
Then $\delta$-invariance of the action
functional and the integration measure guarantees
 that  the addition of a $\delta$-exact term 
to a functional $A$  does not change the 
  correlation function $\Big<A\Big>$. This allows to
identify the space of observables with $\delta$-cohomology.
 Local observables i.e. depending on the
 values of the fields in the vicinity of a given point $(z,\zb)$ can be
described explicitly in terms of the cohomology ring $H^*(X)$. One
has  canonical evaluation map \be ev:\,\,\Sigma\times
\CM(\Sigma,X)\rightarrow X. \ee Then for any $\beta\in H^*(X)$ we
have $ev^*(\beta)\in H^*(\Sigma\times
\CM(\Sigma,X))=H^*(\Sigma)\otimes H^*(\CM(\Sigma,X))$. Let
$ev^*(\beta)=\CO_{\beta}^{(0)}+\CO_{\beta}^{(1)}+\CO_{\beta}^{(2)}$
be a decomposition of $ev^*(\beta)$ with respect to the
grading on the space of  the differential forms  over $\Sigma$.
Then the condition \be\label{decent}
(d+\delta)(\CO_{\beta}^{(0)}+\CO_{\beta}^{(1)}+\CO_{\beta}^{(2)})=0,
\ee implies that the function $\CO_{\beta}^{(0)}$ on $\Sigma$ is
$\delta$-closed and its evaluation at any $p \in \Sigma$  gives 
a cohomology class in $H^*(\CM(\Sigma,X))$. Moreover due to
relation \eqref{decent} we have \be\label{variation}
d\CO_{\beta}^{(0)}=-\delta \CO_{\beta}^{(1)}, \ee and therefore 
this  cohomology class in $H^*(\CM(\Sigma,X))$ does not depend on
the choice of the point $p\in \Sigma$.

For example the local observable at $p\in \Sigma$ constructed using
the  K\"{a}hler form
 \be
  \omega=\sum_{i,\jb=1}^{\ell+1}\,\omega_{i\jb}\,
  d\varphi^j\wedge d\bar{\varphi}^{\jb},
 \ee
on $X$ is given by \be\label{Locobserv}
\CO_{\omega}^{(0)}(p)=\sum_{i,\jb=1}^{\ell+1}\,
\omega_{i\jb}(\varphi(p))\,\chi^i(p)\,\, \bar{\chi}^{\jb}(p). \ee
The corresponding differential form on $\CM(\Sigma,X)$ provides a
K\"{a}hler form on $\CM(\Sigma,X)$.

Correlation functions in topological sigma models are closely
related to  Gromov-Witten invariants counting holomorphic curves
on symplectic manifolds. This can be illustrated as follows.
Taking $t=0$ in \eqref{Gsigma} we obtain the action linearly
dependent  on the fields $F$, $\bar{F}$. Functional integration over
$F$, $\bar{F}$ in \eqref{funcint}  reduces to  the  integration
over the subset of holomorphic maps $\Phi:\Sigma \to X$. Thus the
correlation functions in the topological sigma model at $t=0$ reduce
to counting of holomorphic curves in $X$. Let us notice that
although variation of the parameter $t$  in  \eqref{Gsigma} changes
the action by $\delta$-exact term the correlation functions of
$\delta$-closed functionals depend on $t$ due to holomorphic anomaly
\cite{BCOV}. Thus the interpretation in terms of holomorphic curves
counting in general is not directly applicable for $t\neq 0$.

\subsection{Equivariant  topological linear sigma model}

In the following we need a simple generalization of the standard
notion of the topological sigma-model discussed above. We consider
equivariant version of type A topological linear sigma model on a
disk $D$ with a non-compact target space. In the related
context the mathematical description of correlation functions in
equivariant topological sigma models for compact $X$ and compact
$\Sigma$ was considered
by Givental \cite{Gi1}, \cite{Gi2}, \cite{Gi3}. For  various
application of equivariant cohomology in topological field theory
see e.g. \cite{CMR} and reference therein.

Let us consider  topological sigma model
on the disk $D=\{z|\,\,|z|\leq 1\}$
with   the target space $X=\IC^{\ell+1}$
supplied  with the K\"{a}hler form and K\"{a}hler metric
\be\label{formONE}
\omega=\frac{\imath}{2}\sum_{j=1}^{\ell+1}\,\,d\varphi^j\wedge
d\bar{\varphi}^{j},\qquad
g=\frac{1}{2}\sum_{j=1}^{\ell+1}\,\,(d\varphi^j\otimes
d\bar{\varphi}^{j}+ d\bar{\varphi}^{j}\otimes  d\varphi^{j}).
\ee
We fix the flat metric $h$  on $D$
\be\label{metricD}
h=\frac{1}{2}(dz d\zb+d\zb\,dz)=
(dr)^2+r^2(d\sigma)^2, \qquad r\in [0,1],\quad \sigma\in [0,2\pi],
\ee
where $z=re^{\imath \sigma}$.

The action \eqref{Gsigma} for $X=\IC^{\ell+1}$ and $t=t^*=0$ is
given by 
 \be\label{linsigmact}
  S_D=\int_D\,\,d^2z\,\delta\CV= \imath \int_D\,d^2z\,\,
  \Big(\<F,\apr \varphi\>+ \<\bar{F},\pr
  \bar{\varphi}\>+\<\bar{\psi},\pr\bar{\chi}\>+\<\psi,\apr\chi\>\Big),
 \ee
where $\CV=\imath \<\psi,\apr\varphi\>+\imath
\<\bar{\psi},\pr\bar{\varphi}\>$ and \be\label{transformtwo} \delta
\varphi=\chi,\qquad \delta \chi=0,\qquad \delta \bar{\psi}=
\bar{F},\qquad \delta \bar{F}=0, \ee
$$
\delta \bar{\varphi}=  \bar{\chi},\qquad \delta \bar{\chi}=0,
\qquad \delta \psi= F,\qquad
 \delta F=0.
$$

 To construct an  equivariant extension
of the topological theory with the action \eqref{linsigmact} we
modify  transformations \eqref{transformtwo} following the
interpretation of $\delta$ as a de Rham differential in the
infinite-dimensional setting. Let us first recall the standard
construction of an algebraic model of equivariant cohomology. Let
$M$ be a $2(\ell+1)$-dimensional manifold supplied with an action of
a compact Lie group $G$.  Let $\mathfrak{g}$ be the corresponding
Lie algebra and $\mathfrak{g}^*$ be its dual. Pick a bases $\{t^a\}$,
$a=1,\ldots ,\dim \mathfrak{g}$ in the Lie algebra  $\mathfrak{g}$
and let $\{v^a\}$ be the set of  the corresponding  vector fields on
$M$. Introduce a bases   $\{u^a\}$, $a=1,\ldots ,\dim \mathfrak{g}$
in $\mathfrak{g}^*$ dual to the bases $\{t^a\}$ $a=1,\ldots ,\dim
\mathfrak{g}$. The Cartan  algebraic model of  $G$-equivariant de
Rham cohomology $H^*_G(M)$ is defined as a cohomology of the
following equivariant extension $(\Omega_G^*(M),d_G)$ of the
standard de Rham complex $(\Omega^*(M),d)$ 
\be\label{equivdif}
\Omega^*_G(M)=(\Omega^*(M)\otimes S^*[\mathfrak{g}^*])^G, \qquad
d_G=d-\sum_{a=1}\,u_ai_{v^a}.
\ee
In \eqref{equivdif} the standard 
coadjoint action of $G$ on $\mathfrak{g}^*$ is implied. We have 
\be
d_G^2=-\CL_{\sum_{a=1}^{\dim \mathfrak{g}}\,u_av^a},\qquad
\CL_v=d\,i_v+i_v\,d,
\ee
where $\CL_v$ is the  Lie derivative
along a vector field $v$. Note that  the equivariant
differential $d_G$ satisfies  $d^2_G=0$
 when  acting on $\Omega^*_G(M)$.   The cohomology groups
$H^*_G(M)$ of the complex \eqref{equivdif} naturally have a 
structure of modules  over $H^*_G({\rm pt})=(S^*(\mathfrak{g}^*))^G$
and the algebra $(S^*(\mathfrak{g}^*))^G$ of $G$ invariant
polynomial functions on $\mathfrak{g}$ can be identified with the
algebra
 $(S^*(\mathfrak{h}^*))^W$
of Weyl-invariant functions on  the
 Cartan subalgebra $\mathfrak{h}\subset \mathfrak{g}$.

We have the standard action of  $U_{\ell+1}$ on $V=\IC^{\ell+1}$
and an action of $S^1$ on $D$ by rotations
$\sigma\to \sigma+\alpha$ leaving the metric \eqref{metricD}
invariant. The action of $G=S^1\times U_{\ell+1}$ lifts naturally to the action
on the fields
$(F,\bar{F},\varphi,\bar{\varphi},\psi,\bar{\psi},\chi,\bar{\chi})$.
 For instance the
infinitesimal rotation acts via  the Lie derivative $\CL_{v_0}$
along the vector field $v_0=\frac{\pr}{\pr \sigma}$ 
\be\label{Soneact}
  \delta_{v_0}\varphi^j=\CL_{v_0}\varphi^j=i_{v_0}d\varphi^j,
  \qquad
  \delta_{v_0} F^j=\CL_{v_0} F^j=i_{v_0}dF^j+d(i_{v_0}\,F^j).
 \ee
Let $\hbar v_0$ be an element of the Lie algebra of $S^1$ and
$\Lambda $ be an image of the element $\sum_{a=1}^{\dim
  \mathfrak{g}}u_at^a$ of the Lie algebra $\mathfrak{u}_{\ell+1}={\rm Lie}(U_{\ell+1})$
in the representation $V=\IC^{\ell+1}$.
 The $G$-equivariant analog $\delta_G$ of the
transformation $\delta$  \eqref{transformtwo}
is obtained adapting \eqref{equivdif} to the considered action of
$G=S^1\times U_{\ell+1}$ on the fields of topological sigma model
\be\label{transformeq}
\delta_G \varphi=\chi,\qquad \delta_G
\chi=-(\imath \Lambda\varphi+ \hbar\,\CL_{v_0}\varphi),\qquad
\delta_G \psi=F,\qquad
 \delta_G F=-(\imath\Lambda\psi+\hbar\,\CL_{v_0}\psi),
\ee
\be\nonumber
\delta_G \bar{\varphi}=\bar{\chi},\qquad
\delta_G \bar{\chi}=-(-\imath\Lambda\bar{\varphi}+\hbar\,\CL_{v_0}\bar{\varphi}),
\qquad \delta_G \bar{\psi}=\bar{F},\qquad \delta_G
\bar{F}=-(-\imath \Lambda\bar{\psi}+
\hbar\,\CL_{v_0}\bar{\psi}).
\ee
It is easy to check that the action \eqref{linsigmact} is
both $G$- and  $\delta_G$-invariant.

Equivariant  observables  are given by $\delta_G$-closed
$G$-invariant functionals of the
fields.  The simple direct check shows that the following functional
 defines an equivariant extension of the observable \eqref{Locobserv}
corresponding to  the K\"{a}hler form \eqref{formONE}
\be\label{observone}
\CO_{\Lambda,\hbar}(r)=\frac{\imath}{2}\,\int_0^{2\pi}\,d\sigma \,
(-\<\chi(re^{\imath \sigma}),\chi (re^{\imath \sigma})\>+
\<\varphi(re^{\imath \sigma}), (\imath \Lambda+\hbar \CL_{v_0})
\varphi(re^{\imath \sigma})\>).
\ee

We would like to consider the functional integral over the disk with
the action \eqref{linsigmact} deformed by the observable
$\CO_{\Lambda,\hbar}:=\CO_{\Lambda,\hbar}(r=1)$ where
$\CO_{\Lambda,\hbar}(r)$ is given by \eqref{observone}.  Note that
this functional integral is a Gaussian integral in the
infinite-dimensional space and thus allows mathematically correct
formulation using for instance the $\zeta$-function regularization
\cite{RS}, \cite{H}.

\begin{te} \label{mainTh}
Let $V=\IC^{\ell+1}$ be the  standard representation of
  $U_{\ell+1}$. Let $\Lambda$ be the image of an element $u\in
  \mathfrak{u}_{\ell+1}$ in $End(V)$.  
Then the following identity holds
\be\label{matrixG}
\Big<e^{\mu \,\CO_{\Lambda,\hbar}}\Big>_D=\,\,\hbar^{-\frac{\ell+1}{2}}\,
\det_V\,\,\,\left(\frac{2}{\mu\hbar}\right)^{-\Lambda/\hbar}\,\Gamma(\Lambda/\hbar),
\ee
where $\CO_{\Lambda,\hbar}$ is given by \eqref{observone} for $r=1$.
The  functional integral in the
$S^1\times U_{\ell+1}$-equivariant  Type $A$ topological
linear sigma model with the target space $V=\IC^{\ell+1}$
and the action functional \eqref{linsigmact} is calculated using
$\zeta$-function regularization.
After taking $\mu=2/\pi$, $\hbar=1$  and making the change of
variables $\Lambda\to (s\,\cdot  {\rm id}-\Lambda)/2$ the partition
function \eqref{matrixG} turns into local Archimedean
  $L$-factor \eqref{archLL}.
\end{te}

\noindent{\it Proof}.  The functional integral in \eqref{matrixG}
regularized using $\zeta$-function regularization is
invariant with respect to the action of $U_{\ell+1}$.
Thus the correlation function \eqref{matrixG} is
a central function of $(\ell+1)\times (\ell+1)$ matrix
$\Lambda$ and can be reconstructed from its
restriction to diagonal matrices 
$\Lambda={\rm diag}(\lambda_1,\ldots ,\lambda_{\ell+1})$.
Then the integral is factorized on the product of the
functional integrals in  $S^1\times U_1$-topological sigma models with
one-dimensional target spaces $X=\IC$.
Using the Proposition \ref{intGamma} proved in the next subsection
we have
\be\label{prodG}
\Big<e^{\mu\,\CO_{\Lambda,\hbar}}\Big>_D=
\hbar^{-\frac{\ell+1}{2}}\,
\prod_{j=1}^{\ell+1}\,\left(\frac{2}{\mu\hbar}\right)^{-\lambda_j/\hbar}
\,\Gamma(\lambda_j/\hbar),\qquad
\Lambda=
{\rm diag}(\lambda_1,\ldots \lambda_{\ell+1}).
\ee
Now it is clear that the r.h.s. of \eqref{prodG} is a restriction of
r.h.s. of \eqref{matrixG} to the diagonal matrix $\Lambda=
{\rm diag}(\lambda_1,\ldots \lambda_{\ell+1})$ $\Box$

\begin{rem} According to  Remark \ref{GenLfac},
the correlation function \eqref{matrixG}
for  arbitrary $\mu$ and $\hbar$ can also be  considered as a complex
Archimedean $L$-factor. The $\mu$-dependence of the functional
integral defines also a renormalization scale dependence
(see e.g. \cite{H}). For example using another regularization
of the Gaussian functional integral
(i.e. different form of  $\zeta$-function regularization)
we would obtain the  same  answer up to a  multiplication by
a factor of the form $A\,B^s$. 
\end{rem}

\subsection{Calculation of the functional integral}

In this subsection we prove \eqref{prodG}. It is enough to consider
the case of $\ell=0$. 

\begin{prop} \label{intGamma}
The following integral representation for
  $\Gamma$-function holds
\be\label{InRepG} \Gamma(\lambda/\hbar)=
\hbar^{\frac{1}{2}}\,\left(\frac{2}{\mu\hbar}\right)^{\lambda/\hbar}
\Big<e^{\mu\CO_{\lambda,\hbar}}\Big>_D, \ee where \be
\CO_{\lambda,\hbar}=\,\frac{\imath}{2}\int_{\pr D=S^1}\,d\sigma \,
(-\bar{\chi}\,\chi+
\bar{\varphi}(\imath\lambda+\hbar\pr_{\sigma})\varphi), \ee and the
functional  integral is taken with the action functional 
 \be\label{linsigmactone}
  S_D= \imath\int_{D}\,d^2z\,
  \Big(\bar{F}_z\,\pr_{\zb} \varphi+
  F_{\zb}\,\pr_z \bar{\varphi}+\psi_{\zb}\pr_z\bar{\chi}+\bar{\psi}_z
  \pr_{\zb}\chi\Big).
 \ee
\end{prop}

\noindent{\it Proof}.
Decompose the fields $\varphi$ and $\chi$
\be
\varphi=\varphi_0+\varphi_h,\qquad \chi=\chi_0+\chi_h,
\ee
so that  $\varphi_0$ and $\chi_0$ satisfy the boundary conditions
\be
\varphi_0|_{\pr D=S^1}=0,\qquad \chi_0|_{\pr D=S^1}=0,
\ee
and $\varphi_h$, $\chi_h$ are harmonic functions 
\be
\pr_z\pr_{\zb}\varphi_h=0,\qquad
\pr_z\pr_{\zb}\chi_h=0.
\ee
Let us also decompose the fields $F$ and $\psi$   as follows
$$
\bar{F}_zdz=\bar{G}_zdz+\bar{f}_zdz,
\qquad F_{\zb}d\zb=G_{\zb}d\zb+f_{\zb}d\zb,
$$
$$
\bar{\psi}_zdz=\bar{\xi}_zdz+\bar{\rho}_zdz,
\qquad \psi_{\zb}d\zb=\xi_{\zb}d\zb+\rho_{\zb}d\zb,
$$
where one-forms $f$, $\rho$  and $\bar{f}$, $\bar{\rho}$  satisfy
the equations
\be
\label{fharm} \pr_{\zb}\,\bar{f}_z=0,\qquad \pr_z\,f_{\zb}=0,
\qquad \pr_{\zb}\,\bar{\rho}_z=0,\qquad \pr_z\,\rho_{\zb}=0.
\ee
The fields  $G$, $\bar{G}$ should be in the subspace orthogonal to
the subspace spanned by $f$ and $\bar{f}$ satisfying \eqref{fharm}
\be\label{Forthrel}
\<\bar{G},f\>=\int_D\,d^2z\,\,\bar{G}_z\,f_{\zb}=0,\qquad
\<G,\bar{f}\>=\int_D\,d^2z\,\,G_{\zb}\,\bar{f}_{z}=0 
\ee 
We impose similarly conditions on $\xi$ and $\bar{\xi}$. The following
identity  can be  easily checked \be\label{idetnone}
\int_Dd^2z\,(\bar{F}_z\pr_{\zb}\varphi+F_{\zb}\pr_{z}\bar{\varphi})
=\int_Dd^2z\,(\bar{G}_z\pr_{\zb}\varphi_0+G_{\zb}\pr_{z}\bar{\varphi}_0)
+\int_{\pr D}\,(dz\,\bar{f}_z\varphi_h-d\zb f_{\zb}\bar{\varphi}_h).
\ee
Using \eqref{idetnone} we have the decomposition of the action
\eqref{linsigmactone}
\be\label{decompact}
S_D=S_{bulk}+S_{boundary}, \ee where \be S_{bulk}=\imath
\int_D\,d^2z\,(\bar{G}_z\,\pr_{\zb}\varphi_0+G_{\zb}\,\pr_z\bar{\varphi}_0+
\bar{\xi}_z\,\pr_{\zb}\chi_0+\xi_{\zb}\,\pr_z\bar{\chi}_0), 
\ee
 and
\be S_{boundary}=\imath \int_{S^1=\pr D}\,d\sigma(
\bar{f}_z\,\varphi_h-f_{\zb}\,\bar{\varphi}_h+
\bar{\rho}_z\,\chi_h-\rho_{\zb}\,\bar{\chi}_h). \ee The integration
measure is defined using the standard Hermitian metric on tensor
fields on the disk. For example, the measure
$[D\varphi]\,[D\bar{\varphi}]$  over $(\varphi,\bar{\varphi})$  is
induced by the metric
\be\label{metricF} \|\delta \varphi
\|^2=\int_D d^2z\,\delta\overline{\varphi}\,\,\delta \varphi,\qquad
\|\delta \chi \|^2=\int_D d^2z\,\,\delta\overline{\chi}\,\,\delta
\chi.
\ee
Taking into account that super-manifolds have canonical
integration measure (see  Appendix) the integration measures 
can be split as follows 
 \be
  [DF]\,[D\bar{F}]\,[D\psi]\,[D\bar{\psi}]\,=\,
  [DG]\,[D\bar{G}]\,[Df]\,[D\bar{f}]
  \,[D\xi]\,[D\bar{\xi}]\,[D\rho]\,[D\bar{\rho}],
 \ee
 \be
  [D\varphi]\,[D\bar{\varphi}]\,[D\chi]\,[D\bar{\chi}]
  =\,[D\varphi_0]\,[D\bar{\varphi}_0]\,
  [D\varphi_h]\,[D\bar{\varphi}_h]\,[D\chi_0]\,[D\bar{\chi}_0]\,
  [D\chi_h]\,[D\bar{\chi}_h].
 \ee
Thus using the decomposition \eqref{decompact}
we have a factorization of the functional integral
in \eqref{InRepG} into
 the  product of the integrals over $(G,\varphi_0,\xi,\chi_0)$ and
$(f,\varphi_h,\rho,\chi_h)$. Consider first the integration
over   $(G,\varphi_0,\xi,\chi_0)$ with the action
\be
S_{bulk}=\imath \int\,d^2z\,(\bar{G}_z\,\pr_{\zb}\varphi_0+
G_{\zb}\,\pr_z\bar{\varphi}_0+
\bar{\xi}_z\,\pr_{\zb}\chi_0+\xi_{\zb}\,\pr_z\bar{\chi}_0).
\ee
This is an infinite-dimensional analog of the integral
\eqref{degenerateodd} defined by  $\zeta$-function regularization. The
integral gives a trivial contribution to \eqref{InRepG}.
Thus one should  calculate the following  functional integral
\be
\int [\CD\varphi_h]\,
[\CD\bar{\varphi}_h]\,[\CD
f]\,[\CD\bar{f}]\,[\CD\rho]\,[\CD\bar{\rho}]\,
[\CD\chi_h]\,[\CD\bar{\chi_h}]\,\,e^{-S_*},
\ee
where
\be\label{boundact}
S_*= \int_{S^1=\pr D}\,d\sigma(
\imath\bar{f}_z\,\varphi_h-\imath f_{\zb}\,\bar{\varphi}_h+
\imath\bar{\rho}_z\,\chi_h-\imath\rho_{\zb}\,\bar{\chi}_h
+\mu\frac{\imath}{2}\left(\chi\,\bar{\chi}\,+
\bar{\varphi}(\imath \lambda+\hbar\pr_{\sigma})\varphi\right)).
\ee
Consider decomposition of the fields
\be\label{decone}\nonumber
\varphi_h(z,\zb)=\sum_{n>0}\,\varphi_{-n}\,\zb^n+
\sum_{n\geq 0}\,\varphi_n\,z^n=
\varphi_-(\zb)+\varphi_+(z),
\ee
\be\label{dectwo}
\bar{\varphi}_h(z,\zb)=\sum_{n\geq 0}\,\bar{\varphi}_{-n}\,\zb^n+
\sum_{n> 0}\,\bar{\varphi}_n\,z^n=
\bar{\varphi}_-(z)+\bar{\varphi}_+(\zb),
\ee
\be\label{deconef}\nonumber
\chi_h(z,\zb)=\sum_{n>0}\,\chi_{-n}\,\zb^n+
\sum_{n\geq 0}\,\chi_n\,z^n=
\chi_-(\zb)+\chi_+(z),
\ee
\be\label{dectwof}\nonumber
\bar{\chi}_h(z,\zb)=\sum_{n>0}\,\bar{\chi}_{-n}\,\zb^n+
\sum_{n\geq 0}\,\bar{\chi}_n\,z^n=\bar{\chi}_-(z)+\bar{\chi}_+(\zb).
\ee
Taking into account \eqref{fharm} one can write the action
\eqref{boundact} as follows
\be
S_*= \int_{S^1=\pr D}\,d\sigma(
\imath\bar{f}_z\,\varphi_--\imath f_{\zb}\,\bar{\varphi}_-+
\imath\bar{\rho}_z\,\chi_--\imath\rho_{\zb}\,\bar{\chi}_-
+\mu\frac{\imath}{2}(\chi\,\bar{\chi}\,+
\bar{\varphi}(\imath \lambda+\hbar\pr_{\sigma})\varphi)).
\ee
Integrating over $(\varphi_-,\bar{\varphi}_-,f,\bar{f})$
and $(\chi_-,\bar{\chi}_-,\rho,\bar{\rho})$ and using
\eqref{degenerate} we are left with  the following integral
\be\label{integralLAST}
\Big<e^{\mu\CO_{\lambda,\hbar}}\Big>_D=\int [\CD\varphi_-]\,
[\CD\bar{\varphi}_+]\,[\CD\chi_-]\,[\CD\bar{\chi_+}]\,\,e^{-S_b},
\ee
where
\be\label{boundactLAST}
S_b=- \mu \frac{\imath}{2}\int_{0}^{2\pi}\,d\sigma(\chi_+\,\bar{\chi}_++
\bar{\varphi}_+(\imath \lambda+\hbar\pr_{\sigma})\varphi_+).
\ee
and functional measure is defined using the metric
induced by \eqref{metricF} on the space of the fields \,\,\,\,\,
$(\varphi_+,\bar{\varphi}_+,\chi_+,\bar{\chi}_+)$.

The functional integral \eqref{integralLAST} can be easily
calculated using the $\zeta$-function regularization. Using
the normalization of the integration measure given in Appendix  we
have \be\label{Detexpr}
\Big<e^{\mu\,\CO_{\lambda,\hbar}}\Big>_D=\frac{\det \CD_0}{\det
\,\CD}. \ee Here  the differential operator $\CD=
\,-\frac{\imath}{2}\,\mu\,(\hbar\frac{\pr}{\pr
  \sigma}+\imath \lambda)$  acts in the space of complex-valued functions
on a circle $\sigma\sim   \sigma+2\pi$
which are restrictions of  holomorphic functions on the disk $D$.
The operator $\CD_0$ acts  by  multiplication on $2\pi \cdot
\frac{\mu}{2}$. The spectrum of $\CD$  is
\be
\lambda_n=\frac{\mu}{2}(\hbar n+\lambda),\qquad n\geq 0,
\ee
and using \eqref{regOne}  we have
\be\label{Done}
\ln \det \CD
=(\frac{1}{2}-\frac{\lambda}{\hbar})\ln\left(\frac{\mu \hbar}{2}\right)
+\frac{1}{2}\ln 2\pi-\ln \Gamma(\frac{\lambda}{\hbar}).
\ee
The determinant of $\CD_0$ is calculated  using \eqref{regsum}
\be\label{Dtwo}
\ln \det \CD_0
=\frac{1}{2}\ln\left(2\pi\,\frac{\mu}{2} \right).
\ee
Substitution of \eqref{Done} and \eqref{Dtwo} into \eqref{Detexpr}
gives  \eqref{InRepG} $\Box$

\section{ $\Gamma$-function as an equivariant  symplectic volume}

In the previous Section we represent a product of $\Gamma$-functions
as a  particular Gaussian functional integral. In this section we
interpret this functional integral as an
equivariant symplectic volume of the space of
holomorphic maps of the disk $D$ into $\IC^{\ell+1}$.
According to the general correspondence principle in quantum mechanics
equivariant symplectic volume can be considered as a classical limit
of a partition function of quantum system obtained by a quantization
of the symplectic manifold. We apply this heuristic principle to
the space of holomorphic maps of the disk $D\to \IC^{\ell+1}$. The
corresponding partition function is given by the product of
$q$-deformed $\Gamma$-functions. In the next section we provide a
quantum field theory interpretation of this result by identifying the
product of $q$-deformed $\Gamma$-functions with the correlation
function in the three-dimensional equivariant topological linear sigma-model.

\subsection{$U_{\ell+1}$-equivariant symplectic volume of
  $\IC^{\ell+1}$ }

In this Subsection we consider a
calculation of $U_{\ell+1}$-equivariant symplectic volume of
  $\IC^{\ell+1}$ and its interpretation as
an assymptotic of a partition function of the  associated quantum system. These
considerations  provide  an  example for  the
discussion of  the equivariant symplectic volume of the space of
holomorphic maps of the disk $D$ into $\IC^{\ell+1}$.

Let $M$ be a $2(\ell+1)$-dimensional symplectic manifold
with a symplectic form $\omega$.
Let  $G$ be a compact Lie group acting on $M$.
Let $\mathfrak{g}^*$ be a dual to the Lie
algebra $\mathfrak{g}$ of $G$.  Let 
the action of $G$ on $(M,\omega)$ is Hamiltonian with the momentum map
$H:\,M\to \mathfrak{g}^*$.  
Define $G$-equivariant symplectic volume of
$M$ as follows
\be\label{DH}
Z(M,\lambda)=\int_M\,e^{\<\lambda, H\>+\omega}=
\int_M\,\frac{\omega^{\ell+1}}{(\ell+1)!}\,\,\,e^{\<\lambda, H\>},\qquad
\lambda \in \mathfrak{g},
\ee
where $\<\,,\,\>$ is the  paring between $\mathfrak{g}$ and its dual
$\mathfrak{g}^*$.

The term ``equivariant volume'' for the integral \eqref{DH} comes
from the interpretation of $\omega_G:=\omega+\<\lambda,H\> $ as a 
$G$-equivariant extension of the symplectic form $\omega$.
Indeed, $\omega_G$ is a  $G$-invariant $d_G$-closed
two-form and thus defines an element in $H^2_G(M)$ using  the Cartan model
\eqref{equivdif}  of the equivariant cohomology.  The
interpretation of \eqref{DH} in terms of equivariant cohomology
leads to a possibility to apply a
 powerful localization technique to the calculation the integrals \eqref{DH}
(see \cite{DH} for a direct calculation and
\cite{AB}, \cite{W3} for equivariant localization
approach in abelian and non-abelian cases).

Integral expression \eqref{DH} is invariant with respect to the
adjoint action  of $G$ on $\lambda\in\mathfrak{g}$ and thus can
be uniquely reconstructed from the specialization of \eqref{DH} to
the case when $\lambda$ belongs to a Cartan subalgebra
$\mathfrak{h}\subset \mathfrak{g}$ \be\label{DHCart}
Z(M,\underline{\lambda})=\int_M\,e^{\omega+\sum_{j=1}^{\dim
    \mathfrak{h}}  \lambda_j\, H_j}
,\qquad
\underline{\lambda}=(\lambda_1,\lambda_2,\ldots
,\lambda_{\dim
    \mathfrak{h}}) \in \mathfrak{h}.
\ee

Let us consider a particular example of the equivariant symplectic volume
calculation. Let $M=\IC^{\ell+1}$ be supplied with the symplectic structure
\be\label{symstr}
\omega=\frac{\imath}{2}\sum_{j=1}^{\ell+1}\,dz^j\wedge d\zb^j.
\ee
The symplectic structure \eqref{symstr} is invariant with respect to the
standard action of $U_{\ell+1}$ on $\IC^{\ell+1}$.
The momenta corresponding to the action of the generators
of Cartan subgroup  $(U_1)^{\ell+1}\subset U_{\ell+1}$  are given by
\be\label{Uone}
H_j=-\frac{1}{2}|z^j|^2, \qquad j=1,\ldots ,\ell+1.
\ee
The $(U_1)^{\ell+1}$-equivariant extension
$\omega_{(U_1)^{\ell+1}}$ of $\omega$ satisfies the condition
\be
(d-\sum_{j=1}^{\ell+1}\lambda_ji_{V_j})\omega_{(U_1)^{\ell+1}}=0,
\ee
and is given by
\be
\omega_{U_{\ell+1}}=\omega+\sum_{j=1}^{\ell+1}\lambda_j H_j=
\frac{\imath}{2}\sum_{j=1}^{\ell+1}(
dz^j\wedge d\zb^j+\imath \lambda_j|\varphi^j|^2).
\ee
The simple direct calculation gives
\be\label{fdintegral}
Z(\IC^{\ell+1},\underline{\lambda})
=\int_{\IC^{\ell+1}}\,e^{\omega+\sum_{j=1}^{\ell+1}\lambda_j H_j}
=
\int_{\IC^{\ell+1}}\,\frac{\omega^{\ell+1}}{(\ell+1)!}
\,e^{\sum_{j=1}^{\ell+1}\lambda_j H_j}=
\frac{(2\pi)^{\ell+1}}{\prod_{j=1}^{\ell+1}
\,\lambda_j},
\ee
where $\underline{\lambda}=(\lambda_1,\lambda_2,\ldots
,\lambda_{\ell+1})$ is an element of the Cartan  subalgebra of
$\mathfrak{u}_{\ell+1}={\rm Lie}(U_{\ell+1})$.

It is useful to rewrite the integral \eqref{fdintegral} using
auxiliary  anticommuting variables $(\eta^j$, $\bar{\eta}^{j})$, \\
$j=1,\ldots (\ell+1)$. The standard identification of
polynomial differential forms on $\IC^{\ell+1}$ with  polynomials
in commuting and anticommuting variables
$\IC[z^j,\zb^j,\eta^j,\bar{\eta}^j]$, 
leads to the following expression for 
the symplectic form \eqref{symstr}
 \be\label{symthree}
\omega=\frac{\imath}{2}\sum_{j=1}^{\ell+1}\,\eta^j\, \bar{\eta}^{j}.
\ee The integral \eqref{fdintegral} can be rewritten using the
integration over anticommuting variables in the following form 
 \be\label{fdintegralodd} Z(\IC^{\ell+1},\underline{\lambda})
=\int_{\IC^{\ell+1|\ell+1}}\,\prod_{j=1}^{\ell+1} dz^j d\zb^{j}
\prod_{j=1}^{\ell+1} d\eta^j d\bar{\eta}^{j}\, \,\,\,e^{\frac{\imath
\mu}{2} \sum_{j=1}^{\ell+1}(\eta^j\, \bar{\eta}^{j}+\imath \lambda_j
|z^j|^2)} =\frac{(2\pi)^{\ell+1}}{\prod_{j=1}^{\ell+1} \,\lambda_j},
\ee Here $\IC^{\ell+1|\ell+1}$ shall be considered as an odd tangent
bundle $\Pi T\IC^{\ell+1}$ to $\IC^{\ell+1}$. Note that the integral
is independent on $\mu$.

According to the Correspondence Principle in quantum
mechanics/statistical mechanics the integral over symplectic
manifolds of the form \eqref{DH} describes an asymptotic  of
the partition function of the corresponding quantum system. More
precisely, if a quantization of the symplectic manifold $(M,\omega)$
exists  and $\hat{H}_j$ are  quantum counterparts of the
Hamiltonians $H_j$ then the equivariant volume \eqref{DHCart}
describes the asymptotic of the following trace over the
Hilbert space of the corresponding quantum system
 \be\label{charep}
Z_{\tilde{\hbar}}(M,\beta,\underline{\lambda})=\Tr_{\CH}\,\,\,\,e^{-\frac{
    \beta}{\tilde{\hbar}}
 \sum_{j=1}^{\ell+1}\lambda_j \hat{H}_j},
\ee
where $\tilde{\hbar}$ is parameter of quantum deformation.

For example in the case of $M=\IC^{\ell+1}$ supplied  with the
symplectic form \eqref{symstr} the corresponding quantum system can
 be easily found. The quantum counterpart of the algebra of polynomial
functions $\IC[z^1,\ldots, z^{\ell+1},\zb^1,\ldots ,\zb^{\ell+1}]$ is
a Heisenberg algebra generated by
$\hat{z}^1,\ldots, \hat{z}^{\ell+1},\hat{\zb}^1,\ldots
,\hat{\zb}^{\ell+1}$ with the relations
\be
[\hat{z}^i,\hat{z}^j]=0,\qquad
 [\hat{\zb}^i,\hat{\zb}^j]=0,\qquad [\hat{\zb}^i,\hat{z}^j]=
-2\tilde{\hbar} \delta^{ij}.
\ee
This algebra can be represented in the space of polynomials
$\CV^{(\ell+1)}=\IC[z^1,\ldots, z^{\ell+1}]$ with the action of the generators
$$
\hat{z}^j=z^j,\qquad \hat{\zb}^j=-2\tilde{\hbar} \frac{\pr}{\pr z^j}.
$$
The space  $\CV^{(\ell+1)}$ has a natural structure of
$U_{\ell+1}$-module  and $U_{\ell+1}$-character of the module $\CV^{(\ell+1)}$ is
easy to calculate. The generators of the Cartan subalgebra of
$\mathfrak{u}_{\ell+1}$ are realized by the following differential
operators \be\label{quantumHam} \hat{H}_j=\tilde{\hbar}
z^j\frac{\pr}{\pr z^j},\qquad j=1,\ldots ,\ell+1, \ee and provide
a  quantization of the classical Hamiltonians \eqref{Uone}. The
character of $U_{\ell+1}$-module $\CV^{(\ell+1)}$ is given by
\eqref{charep} and a simple calculation gives
\be\label{charepone}
Z_{\tilde{\hbar}}(\IC^{\ell+1},\beta,\underline{\lambda})=\Tr_{\CV^{\ell+1}}\,e^{
  -\frac{\beta}{\tilde{\hbar}}\,\sum_{j=1}^{\ell+1}\lambda_j \hat{H}_j}=
\prod_{j=1}^{\ell+1}\,\,
\frac{1}{1-e^{-\beta \lambda_j}}.
\ee
Then the  following relation between equivariant volume
\eqref{fdintegral} and the  character \eqref{charepone} holds
\be\label{assymptFD}
Z(\IC^{\ell+1},\underline{\lambda})=\lim_{\beta \to
  0}\,(2\pi \beta)^{\ell+1}\,Z_{\tilde{\hbar}}(\IC^{\ell+1},\beta,\underline{\lambda}).
\ee

\subsection{$S^1\times U_{\ell+1}$-equivariant volume of
  $\CM(D,\IC^{\ell+1})$}

In this Subsection we provide an interpretation of the
functional integral calculation in Section 2
as a calculation of an equivariant volume
of a space of holomorphic maps of the disk $D$ into $\IC^{\ell+1}$.

Let $\CM(D,\IC^{\ell+1})$ be the  space
of holomorphic map of $D=\{z|\,\,|z|\leq 1\}$ into $\IC^{\ell+1}$.
An element of $\CM(D,\IC^{\ell+1})$
can be described  by a set of functions
$(\varphi^j(z,\zb),\bar{\varphi}^{j}(z,\zb))$,  satisfying the constraints
\be
\pr_{\zb}\varphi^j(z,\zb)=0,\qquad
\pr_{z}\bar{\varphi}^{j}(z,\zb)=0,
\qquad j=1,\ldots
(\ell+1).
\ee
Define a  symplectic form on the space $\CM(D,\IC^{\ell+1})$ of
holomorphic maps  as follows
\be
\Omega=\frac{\imath}{2}
\sum_{j=1}^{\ell+1}\,\int_0^{2\pi}\,\delta \varphi^{j}(\sigma)\wedge
\delta\bar{\varphi}^{j}(\sigma)\,d\sigma,
\ee
where $\varphi^j(\sigma)$, $\bar{\varphi}^{j}(\sigma)$
are restrictions of $\varphi^j(z,\zb)$, $\bar{\varphi}^{j}(z,\zb)$
to  the boundary $\pr D=S^1$.
This symplectic form is invariant with respect to the group $S^1$ of loop rotations
$\sigma\to \sigma +\alpha$ and the  action of $U_{\ell+1}$ is induced from the
standard action on $\IC^{\ell+1}$. The action of $S^1\times
U_{\ell+1}$ is Hamiltonian and the
momenta corresponding to the action of $S^1$ and of the  Cartan
subgroup $U_1^{\ell+1}\subset U_{\ell+1}$
are  given by
\be\label{Momenta}
H_0=\frac{\imath}{2}\sum_{j=1}^{\ell+1}\int_0^{2\pi}\, \bar{\varphi}^{j}
(\sigma)\pr_\sigma \varphi^{j}(\sigma)\,
\,d\sigma,\qquad
H_j=-\frac{1}{2}
\int_0^{2\pi}\, |\varphi^j(\sigma)|^2\,\,d\sigma, \qquad
j=1,\ldots \ell+1.
\ee
In the  infinite-dimensional setting  we can try to define an
analog of the equivariant volume \eqref{DHCart}. We use an analog of
the representation \eqref{fdintegralodd} in terms of the integral over
commuting and anticommuting variables.
Let $\chi^j(z,\zb)$ and $\bar{\chi}^{j}(z,\zb)$ be anticommuting
counterparts  of $\varphi^j(z,\zb)$
$\bar{\varphi}^{j}(z,\zb)$ satisfying the equations
\be
\pr_{\zb}\chi^j(z,\zb)=0,\qquad
\pr_{z}\bar{\chi}^{j}(z,\zb)=0,
\qquad j=1,\ldots
(\ell+1).
\ee
Let $\hbar$ be a
generator of Lie algebra of $S^1$ and $\underline{\lambda}=(\lambda_1,\ldots
\lambda_{\ell+1})$ be generators of Cartan subalgebra of $\mathfrak{u}_{\ell+1}$.
The following formal integral should be considered as a $S^1\times
U_{\ell+1}$-equivariant volume of the space of holomorphic maps
$\CM(D,\IC^{\ell+1})$
\be\label{intSI}
Z(\CM(D,\mu,\IC^{\ell+1}),\beta ,\underline{\lambda})=
\int_{\Pi T\CM(D,\IC^{\ell+1})}\,dm(\varphi,\chi)\,\,\,
\,e^{\frac{\imath\mu}{2}\sum_{j=1}^{\ell+1}\,(\int_0^{2\pi}\,d\sigma\,\chi^j\,
\bar{\chi}^{j}\,+\imath \lambda_j H_j)+\hbar\,H_0},
\ee
where $H_0$, $H_j$ are given by \eqref{Momenta} and
$dm(\varphi,\chi)$ is an integration measure
 to be defined.
This functional integral is a product of
functional integrals \eqref{integralLAST}, \eqref{boundactLAST}
arising at an intermediate step of the proof of  Proposition
\ref{intGamma}. Thus we can define the measure
$d m(\varphi,\chi)$ in \eqref{intSI} in the same way as in
Proposition \ref{intGamma}. Finally taking  into account Theorem \ref{mainTh}
 we arrive to the conclusion that the r.h.s. of the relation
 \eqref{matrixG} is indeed can be considered as  $S^1\times
 U_{\ell+1}$-equivariant volume of the  space of holomorphic maps
$\CM(D,\IC^{\ell+1})$. Thus taking into account the definition of
the equivariant symplectic volume given above, 
 Proposition 2.1 can be reformulated as
follows. {\it  The classical $\Gamma$-function is equal to the 
$S^1\times U(1)$-equivariant  symplectic volume of the space 
$\CM(D,\IC)$ of holomorphic maps of the disk $D$ to the complex plane
$\IC$}. 

Let us push further the analogy with the finite-dimensional case
considered in the previous subsection. We would like to find a
quantum system such that the equivariant volume \eqref{intSI}
appears in the assymptotic of  its  partition function.
This is easy to accomplish taking into account that underlying
infinite-dimensional symplectic  space $\CM(D,\IC^{\ell+1})$ is a
linear space. Consider the Fourier series expansion
\be
\varphi^j(\sigma)=\sum_{n\geq 0 }\,\,(2\pi)^{-1/2}\varphi^j_n\,e^{\imath n\sigma},
\qquad \bar{\varphi}^{j}(\sigma)=\sum_{n\geq
  0}\,(2\pi)^{-1/2}\,\bar{\varphi}^{j}_n\,e^{-\imath n\sigma}.
\ee In terms of the  Fourier coefficients we have
\be\label{sympstrF}
\Omega=\frac{\imath}{2}\sum_{j=1}^{\ell+1}\sum_{n\geq 0}\,\delta
\varphi^{j}_n\wedge \delta \bar{\varphi}_n^{j}, \ee
\be\label{Hamiltonians} H_0=-\frac{1}{2}\sum_{j=1}^{\ell+1}\sum_{n
>0}^{\infty}\,n|\varphi^j_{n}|^2,\qquad H_j=-\frac{1}{2}\sum_{n\geq
0}^{\infty}\,|\varphi^j_{n}|^2. \ee Quantization  of the algebra of
polynomial functions on $\CM(D,\IC^{\ell+1})$ with the symplectic
structure \eqref{sympstrF} is given by \be
[\varphi^j_n,\varphi^k_m]=0,\qquad [\bar{\varphi}^j_n,\varphi^k_m]=
-2\tilde{\hbar}\,\delta^{jk}\delta_{nm}, \qquad
[\bar{\varphi}^j_n,\bar{\varphi}^k_m]=0. \ee 
We chose a polarization 
in such a  way that the Hilbert space $\CH^{(\ell+1)}$ is  realized
as a space of functions of  $\varphi^j_n$, $n \geq 0$, $j=1,\ldots
,(\ell+1)$. Then $\bar{\varphi}^{j}_n$, $n \geq 0$, $j=1,\ldots
,(\ell+1)$  act as first order differential operators \be
\bar{\varphi}^{j}_n=-2\tilde{\hbar}\,\frac{\pr}{\pr \varphi^j_n}.
\ee
 The Hilbert space $\CH^{(\ell+1)}$ has a natural action of $S^1\times
(U_1)^{\ell+1}$ induced by the action on the linear coordinates
$$
e^{\imath \a_0}:\,\,\varphi_n^j\rightarrow e^{\imath
  n\a_0}\varphi^j_n,\qquad e^{\imath\a_0}\in S^1,
$$
$$
e^{\imath \alpha_k}:\,\,\varphi_n^j\rightarrow e^{\imath
  \delta_{k j}\alpha_j}\varphi^j_n\,\qquad e^{\imath\alpha_k}\in U_1.
$$
The action of the corresponding Lie algebra is generated by the vector
fields
\be\label{generators}
\hat{H}_0=\tilde{\hbar}
\sum_{j=1}^{\ell+1}\sum_{n=0}^{\infty}\,n\varphi^j_{n}\frac{\pr}{\pr
  \varphi_n^j},\qquad
\hat{H}_j=\tilde{\hbar}\sum_{n=0}^{\infty}\,\varphi^j_{n}\,\frac{\pr}{\pr
  \varphi_n^j}.
\ee
Thus the infinite-dimensional counterpart of \eqref{charepone}
 is given by  the following partition function
\be\label{partfunc}
Z_{\tilde{\hbar}}(\CM(D,\IC^{\ell+1},\beta,\underline{\lambda})=
\Tr_{\CH^{(\ell+1)}}\,e^{-\frac{\beta}{\tilde{\hbar}}
(\hbar\hat{H}_0+\sum_{j=1}^{\ell+1}\lambda_j\hat{H}_j)},
\ee
where $\underline{\lambda}=(\lambda_1,\ldots ,\lambda_{\ell+1})$.

\begin{prop}\label{Lemma3d} The following identity holds
\be
Z_{\tilde{\hbar}}(\CM(D,\IC^{\ell+1},\beta,\underline{\lambda})
=\prod_{j=1}^{\ell+1}\,\Gamma_q(t_j),
\ee
where $q=e^{-\beta \hbar}$, $t_j=e^{-\beta \lambda_j}$
and the $q$-version of $\Gamma$-function is given by
\be
\Gamma_q(t)=\prod_{k=0}^{+\infty}\,\frac{1}{1-tq^k}.
\ee
 \end{prop}

\noindent{\it Proof}. Direct calculation $\Box$

Note that the first order differential operators \eqref{generators}
can be considered as a quantization of the Hamiltonians
\eqref{Hamiltonians}. Moreover,  similarly to \eqref{assymptFD} the
asymptotic  of \eqref{partfunc} is proportional to the
equivariant symplectic volume \eqref{intSI}. Thus the proposed 
$q$-version of $\Gamma$-function is obtained by applying the
standard quantization rules to the classical $\Gamma$-function
expressed as an equivariant volume integral.

\begin{rem} The Hilbert space $\CH^{(\ell+1)}$ can be considered as
a Hilbert space of a second quantized system i.e.
 $\CH^{(\ell+1)}=S^*V$ is a space of multi-particle states
where $V$ is a one-particle Hilbert space. Explicitly
$\CH^{(\ell+1)}=\IC[\varphi_n^j]$, $n\geq 0$, $j=1,\ldots (\ell+1)$
is  a space of polynomial functions on $\IC^{\ell+1}\otimes \IC[\xi]$ and
$V \subset  \CH^{(\ell+1)}$ is a subspace of linear functions on
$\IC^{\ell+1}\otimes \IC[\xi]$.
\end{rem}

\section{ $q$-deformed $\Gamma$-function:  Three-dimensional interpretation}

In the previous Sections we propose the integral
 representation of the $\Gamma$-function
as a functional integral in the equivariant topological
two-dimensional linear sigma model.  This integral representation
can be interpreted as an equivariant symplectic volume of the space
of holomorphic maps of a disk into a complex plain $\IC$. We also
define a natural quantum version of the equivariant symplectic
volume. This is equal to a $q$-deformed $\Gamma$-function expressed
as a partition function of a quantum system. It is natural to
look for a topological field theory expression for $q$-deformed
$\Gamma$-function similar to the expression for classical
$\Gamma$-function given in Section 2. In this Section we propose
such a formulation in terms of an  equivariant
topological three-dimensional linear sigma model. This
is an example of the representation of a $K$-theory analog of Gromov-Witten
invariants \cite{Gi3}, \cite{GiL} in terms of 
 three-dimensional topological field theories discussed below.

\subsection{Path integral representation for the character}

We start with a path integral representation
of the partition function \eqref{charepone}
of the quantum system  obtained by a quantization
of the classical phase space $\IC^{\ell+1}$ supplied with the symplectic form
\eqref{symstr}. We have the following standard
representation
\be\label{standrrep}
Z_{\tilde{\hbar}}(\IC^{\ell+1},\beta,\underline{\lambda})=
 \int_{L\Pi T(\IC^{\ell+1})}\,[Dz]\,[D\zb] [D\eta]\,[D\bar{\eta}]\,\,
e^{S_{QM}(z,\zb,\eta,\bar{\eta})} , \ee where \be\label{QMaction}
S_{QM}=\frac{1}{\tilde{2\hbar}}\,\,\sum_{j=1}^{\ell+1}\,\,
\int_{0}^{2\pi}\,d\tau \,\left(\zb^j\frac{dz^j}{d\tau}
+\frac{\imath\beta}{2\pi}(\eta_j\bar{\eta}_j+\imath
\lambda_j|z^j|^2)\right). \ee Here commuting periodic functions
$(z^j(\tau),\zb^j(\tau))$, $j=1,\ldots ,(\ell+1)$ and anticommuting
periodic functions  $(\eta^j(\tau), \bar{\eta}^j(\tau))$,
$j=1,\ldots ,(\ell+1)$ parametrize the  loop space of the odd
tangent bundle $\Pi T\IC^{\ell+1}=\IC^{\ell+1|\ell+1}$ to
$\IC^{\ell+1}$. The role of the integral over anticommuting
variables is to provide a proper  integration measure over commuting
variables (see Appendix).

The path integral \eqref{standrrep} is a Gaussian one. By
the definition we use,  it is expressed through the determinants of the
differential operators regularized using the corresponding $\zeta$-functions.
The integral \eqref{standrrep} is easily calculated
using the expression \eqref{regTwo}  for regularized determinants
\be
Z_{\tilde{\hbar}}(\IC^{\ell+1},\underline{\lambda})=\prod_{j=1}^{\ell+1}\,\,
\frac{1}{\det( \frac{d}{d\tau}-\beta\lambda_j/2\pi)}
=\prod_{j=1}^{\ell+1}\,\,
\frac{1}{1-e^{- \beta\lambda_j}}.
\ee
This coincides with the result of the direct calculation
 of the partition function \eqref{charepone}.
Let us remark that the  $\zeta$-function regularization 
implies the following normal ordering: given a monomial in 
$z^j,\bar{z}^j,\,j=1,\ldots,\ell+1$, one should put $\bar{z}^j$'s to 
the l.h.s. and $z^j$'s to the r.h.s.

The functional integral \eqref{standrrep} can be  related with
equivariant cohomology of the loop space $L\IC^{\ell+1}$ (see e.g.
\cite{CMR}). Let us  identify differential forms on
$L\IC^{\ell+1}$ with functions on the loop space of
$\IC^{\ell+1|\ell+1}$. Let $(z^j(\tau),\zb^j(\tau))$, $j=1,\ldots
,(\ell+1)$  and anticommuting functions  $(\eta^j(\tau),
\bar{\eta}^j(\tau))$, $j=1,\ldots ,(\ell+1)$ be linear coordinates
on $L(\Pi T\IC^{\ell+1})$. The standard de Rham differential  on
$L\IC^{\ell+1}$ is then given by \be \delta
z^j(\tau)=\eta^j(\tau),\qquad \delta\eta^j(\tau)=0, \qquad \delta
\zb^j(\tau)=\bar{\eta}^j(\tau),\qquad \delta\bar{\eta}^j(\tau)=0.
\ee
 We have a natural diagonal action of $(U_1)^{\ell+1}$ on $L(\Pi
T\IC^{\ell+1})$ induced by the diagonal action on $\IC^{\ell+1}$.
Let $S^1_{\tau}$ act by rotations of the loops $\tau\to \tau
+\alpha$, $\alpha \in [0,2\pi]$ and let $\tilde{\beta}$ be a
generator of the corresponding Lie algebra. Let
$G_0=S^1_{\tau}\times (U_1)^{\ell+1}$ be an  abelian subgroup of
$G=S^1_{\tau}\times U_{\ell+1}$. Then $G_0$-equivariant de Rham
differential is given by \be\nonumber \delta_{G_0}
z^j(\tau)=\eta^j(\tau),\qquad \delta_{G_0}\eta^j(\tau)=-\left(
\tilde{\beta}^{-1}\frac{dz^j}{d\tau}+\imath \lambda_jz^j\right), \ee
\be\nonumber \delta_{G_0} \zb^j(\tau)=\bar{\eta}^j(\tau),\qquad
\delta_{G_0}\bar{\eta}^j(\tau)=
-\left(\tilde{\beta}^{-1}\frac{d\bar{z}^j}{d\tau}-\imath
\lambda_j\zb^j\right), \ee and satisfy the relation \be
\delta_{G_0}^2=-\tilde{\beta}^{-1}\pr_{\tau}. \ee 
The action functional \eqref{QMaction} is proportional to the $G_0$-equivariant
extension $\omega_G$ of the symplectic form \eqref{symthree} \be
\omega_{G_0}=\frac{\imath }{2}\sum_{j=1}^{\ell+1}\int_0^{2\pi}\,
d\tau\,\left(\bar{\eta}^j\,\eta^j+\imath \lambda_j|z^j|^2+
\tilde{\beta}^{-1}\zb^{j} \frac{d z^j}{d\tau}\right), \ee where
$\tilde{\beta}=\imath \beta/2\pi$. Thus the functional integral with
the action \eqref{QMaction} can be considered as  (an analytic
continuation of) a 
 $G_0$-equivariant symplectic volume of the loop space $L\IC^{\ell+1}$.

\subsection{Three-dimensional interpretation of $K$-theory
Gromov-Witten invariants}

Givental in \cite{Gi3} (see also \cite{GiL}) proposed a construction
of quantum K-theory invariants of complex manifolds
generalizing quantum cohomology invariants arising in Gromov-Witten
theory.  The Givental construction of quantum K-theory invariants
of  a manifold $X$ is formulated directly in terms
of the characteristic classes of complexes of coherent sheaves on
compactified moduli space of holomorphic curves in $X$. For example
let $X$ be a Fano manifold with the cohomology ring generated by
the second cohomology $H^2(X)$, $\dim H^2(X)=m$. Let
$\CM_{\underline{d}}(\Sigma,X)$
be a compactified moduli space of holomorphic maps
of an algebraic curve
$\Sigma$ into $X$ of a fixed multi-degree
$\underline{d}=(d_1,\cdots ,d_{m})\in \IZ^{m}$.
The corresponding invariant is  a generating function
of holomorphic  Euler characteristics  of line bundles on
$\CM_{\underline{d}}(\Sigma,X)$
\be\label{Kinvar}
\CG(X,\underline{Q},\CL)=\sum_{\underline{d}\in \IZ^m}\,
\,{\rm ch}(\CM_{\underline{d}}(\Sigma,X),\CL)\,\prod_{a=1}^m Q_a^{d_a},
\ee
where
$$
{\rm ch}(\CM_{\underline{d}}(\Sigma,X),\CL)=\sum_n\,(-1)^n\,\dim
H^n(\CM_{\underline{d}}(\Sigma,X),\CL),
$$
and $\CL$ is a line bundle on  $\CM_{\underline{d}}(\Sigma,X)$.

In this Subsection we propose a topological field theory
interpretation of the quantum $K$-theory invariants. Les us start
with some heuristic arguments. Cohomology of a manifold $X$ can be
described in  terms of differential forms on $X$ via the de Rham complex. 
On the other hand $K$-cohomology of $X$ can be 
described in terms of vectors bundles on $X$. Thus a
relation between cohomology and $K$-theory  can  be considered as a
kind of a quantization. The general approach to  quantization is
via Feynman path integral. The Feynman path integral  description of quantum
mechanics implies that a quantization of a symplectic manifold
should be described in terms of the geometry of the loop
space on the symplectic manifold. Therefore $K$-theory  of $X$ (modulo
torsion) should  allow interpretation in terms of cohomology of the
loop space $LX$.  Let us apply this reasoning  to quantum
counterpart of cohomology and $K$-theory. 
The quantum cohomology of a manifold $X$ are  
described in terms of a topological sigma model in two dimensions
 with the target space $X$ \cite{W2}. Thus it would be natural to  guess that quantum
$K$-theory invariants should be described in terms of 
two-dimensional sigma model on the loop space $LX$ i.e 
in terms of  three-dimensional topological sigma models. We
propose such  three-dimensional topological theory description below.

Let $\Sigma$ be a Riemann surface and $X$ be a K\"{a}hler manifold.
We retain the  basic notations  introduced in Section 2.1. Let us 
 consider a  three-dimensional manifold  $\Sigma\times
S^1$. We pick  a metric
$d^2s_{\Sigma}=h_{z\zb}\,dz\,d\zb$ on $\Sigma$ and a constant metric
$d^2s_{S^1}=d\tau^2$ on $S^1$. Let $\pi:\Sigma\times S^1 \to \Sigma$
be a projection. Introduce the fields $F$ and $\bar{F}$ as sections of
$\pi^*(K)\otimes \Phi^*(T^{0,1})$ and of $\pi^*(\bar{K})\otimes
\Phi^*(T^{1,0})$ correspondingly. The anticommuting  fields  $\chi$,
$\bar{\chi}$ are sections  of the bundles $\Phi^*(\Pi T^{1,0}X)$,
$\Phi^*(\Pi T^{0,1})$ and anticommuting fields  $\psi$, $\bar{\psi}$
are sections of the bundles $\pi^*(K)\otimes \Phi^*(\Pi T^{0,1})$,
$\pi^*(\bar{K})\otimes \Phi^*(\Pi T^{1,0})$ respectively.  Metrics
$g$ on $X$ and $h$ on $\Sigma$ induce a Hermitian paring $\<\,,\,\>$
on the space of sections  of the considered bundles. An example of
the pairing  in local coordinates is given by \eqref{metrex}.
Consider the following action functional in three dimensions
\be\label{Gaction3d} S^{3d}_0(t)=\beta\,\int_{S^1\times
\Sigma}d\tau\,d^2z\sqrt{h}\, \Big(\imath \<\CF,\apr
\varphi\>+\imath\<\bar{\CF},\pr \bar{\varphi}\> +\imath
\<\bar{\psi}, D\bar{\chi}\>+\imath\<\psi,\bar{D}\chi\>
+t\<\CF\,\CF\> \ee \be\nonumber +\frac{t}{2}\<\psi,\pr_{\tau}\psi\>+
\frac{t}{2}\<\psi,\pr_{\tau}\psi\>+ t\<\bar{\psi},
R(\psi,\bar{\chi})\chi\>+
2\tilde{t}\beta^{-2}\<\pr_{\tau}\varphi,\pr_{\tau}\varphi\>+\imath
\tilde{t}\beta^{-1} \<\chi,\pr_{\tau}\chi\>+\imath
\tilde{t}\beta^{-1} \<\bar{\chi},\pr_{\tau}\bar{\chi}\>\Big), \ee
The fields $\CF$, $\bar{\CF}$  and the derivatives $D_z$,
$\bar{D}_{\zb}$ are locally defined as in \eqref{defone},
\eqref{deftwo}. The action functional \eqref{Gaction3d} is invariant
with respect to the transformations \be\label{transform} \delta
\varphi=\chi,\qquad \delta \chi=-\imath \beta^{-1}
\pr_{\tau}\varphi,\qquad \delta \bar{\psi}= \bar{F},\qquad \delta
\bar{F}=-\imath \beta^{-1} \pr_{\tau}\bar{\psi}, \ee
$$
\delta \bar{\varphi}=  \bar{\chi},\qquad \delta \bar{\chi}=-\imath \beta^{-1} \pr_{\tau}
\bar{\varphi}, \qquad \delta \psi= F,\qquad
 \delta F=-\imath \beta^{-1} \pr_{\tau}\psi,
$$
and the following relation holds
\be\label{equivar3d}
 \delta^2=-\imath \beta^{-1}\pr_{\tau}.
\ee
Similarly to the action $S_0(t)$ in \eqref{Gsigma}
the three-dimensional action \eqref{Gaction3d}
can be written as follows
\be\label{Qcom}
S^{3d}_0(t)=\beta \int_{S^1\times \Sigma}d\tau\,d^2z\sqrt{h}\,\delta\CV(t)
\ee
where
$$\CV(t)=\<\psi,\frac{1}{2}t\CF+\imath \apr\varphi\>+
\imath \tilde{t}\beta^{-1}\<\chi,\pr_{\tau}\varphi\>+
\<\bar{\psi},\frac{1}{2}t\bar{\CF}+\imath \pr\bar{\varphi}\>
+\imath \tilde{t}\beta^{-1}\<\bar{\chi},\pr_{\tau}\bar{\varphi}\>.
$$
Note that the form of the $\delta$-transformations allows to interpret
the three-dimensional sigma model  \eqref{Gaction3d}
in terms of  $S^1$-equivariant cohomology  where $S^1$ acts by shits
of $\tau$.

\begin{rem}
The action \eqref{Gaction3d}  can be obtained by a twisting of 
a three-dimensional  $\CN=2$ SUSY sigma model on
$\Sigma\times S^1$ with the target space $X$ similarly to the
two-dimensional case \cite{W2}. 
\end{rem}

We consider the theory with the action \eqref{Gaction3d} deformed by
a  $\delta$-closed functional.   Let $\CL$ be a holomorphic line
bundle on $X$ with a curvature two-form $\omega^{\CL}$ representing
the first Chern class $c_1(\CL)$. Locally we can express the
curvature  in terms of a connection one-form $\alpha$
\be\label{locsol} \omega^{\CL}=\apr \alpha. \ee The following
functionals are invariant with respect to transformations
\eqref{transform} \be\label{observthree} \exp(2\pi \imath
\CO(\CL))=\exp\left(2\pi \imath \int_{\Sigma \times
S^1}\,\,\,d\tau\, \sqrt{h}\,d^2z\,\Big(\sum_{i,\jb=1}^{\ell+1}
\omega^{\CL}_{i\jb}\chi^i\chi^{\jb}+\sum_{i=1}^{\ell+1}\alpha_i(\varphi)
\pr_{\tau}\varphi^i\Big)\right), \ee \be\label{observablefour}
\exp(\CO_{\omega}(y))=\exp(\sum_{a=1}^{\rm{rk}\,\, H^2(X)}\,y_a
\int_0^{2\pi} d\tau \int_{\Sigma}\,\Phi^*(\omega_a)). \ee where
$\{\omega_a\}$ is a bases  in $H^2(X)$. Note that the
observable \eqref{observthree} is well-defined. Indeed, although it
is written using the local representation \eqref{locsol} the
ambiguity is given by the exponent of the $2\pi \imath $ multiplied
by a period of the two-form $\omega^{\CL}\in H^2(X,\IZ)$  and thus
irrelevant.

The three-dimensional sigma model with the action \eqref{Gaction3d} 
is not well-defined in general as a quantum field theory. The sigma
model provides a long wave length (``infra-red'') 
description and needs a proper short wave length (``ultra-violent'')  completion  in
general. For the case of the target space being a flag manifold
(relevant to the quantum $K$-theory invariants considered in \cite{GiL}) 
the proper completion can be given for example in terms of 
quiver gauge theories with $\CN=4$ SUSY. 
We leave the detailed construction of these completions to
another occasion. Here we provide heuristic arguments to support the
conjecture that the topological field theory with the action
\eqref{Gaction3d} specialized to $t=0$ calculates quantum K-theory
invariants \eqref{Kinvar}.
The argument is very close to that for the standard description of
Gromov-Witten invariants
in terms of a two-dimensional twisted $\CN=2$ SUSY sigma model
\cite{W2}, \cite{AM}. Integrating over $F$ and $\psi$
 we reduce  the integration over other fields
to the vicinity of the subspace of maps $\Sigma\times  S^1\to X$
holomorphic along $\Sigma$
\be\label{holsub}
\pr_{\zb}\varphi^i=0,\qquad D_{\zb}\chi^i=0.
\ee
The corresponding
determinant contributions arising from the integration over
commuting and anticommuting  fields in the quadratic approximation
around the subspace \eqref{holsub}  cancel each other (we consider
the situation when there is no zero modes of $F$ and $\psi$).   Thus
the integral reduces to the path integral with the  following action
\be\nonumber
  S_{red}(\CL)=\,\int_0^{2\pi}\,d\tau\,\left[
  \int_{\Sigma}\sqrt{h} \,\,d^2z\left(\sum_{i,\jb=1}^{\ell+1}
  \,\tilde{t}\,\,\omega^{\CL}_{i\jb}\,\,\,\Big(\beta^{-2}\frac{\pr\,\varphi^i}{\pr
  \tau} \frac{\pr\,\bar{\varphi}^{\jb}}{\pr \tau}
  +\beta^{-1}\bar{\chi}^{\jb}\frac{\pr \chi^{i}}{\pr
  \tau}\Big)\right.\right.
 \ee
 \be\label{reducedAct}
  \left.\left.+\Big(\sum_{i,\jb=1}^{\ell+1}
  \omega^{\CL}_{i\jb}\chi^i\bar{\chi}^{\jb}+
  \sum_{i=1}^{\ell+1}\alpha_i(\varphi) \pr_{\tau}\varphi^i\Big)
  \right)\,+\,
  \sum_{a=1}^{\rm{rk}\,\, H^2(X)}\,y_a
  \int_{\Sigma}\,\Phi^*(\omega_a)\right]\,,
 \ee where $\varphi:\Sigma \to
X$ is a holomorphic map and $\chi$ is a section of the odd tangent
bundle $\Pi T\CM(\Sigma,X)$ to the space $\CM(\Sigma,X)$ of the
holomorphic maps. The theory described by this action is a $\CN=1/2$
SUSY one-dimensional sigma model with the target space
$\CM(\Sigma,X)$  (first line in \eqref{reducedAct}) deformed by  an
observable (second line in \eqref{reducedAct}).  Using the standard
facts on the partition functions for $\CN=1/2$ SUSY quantum
mechanics (see \cite{A}, \cite{FW}, \cite{ASW}, \cite{B})
 we can identify  the result of taking functional integral
with the action \eqref{reducedAct} with the
generating function \eqref{Kinvar} given by the
sum of  holomorphic Euler characteristics of the line bundle $\CL$
with $\omega=c_1(\CL)$  over components $\CM_{\underline{d}}(\Sigma,X)$ of the
moduli space of holomorphic maps with the identification $Q_j=\exp y_j$.
It is easy to generalize these considerations to the case of
an  equivariant  three-dimensional sigma-model taking into account
proper ultra-violent completion. 
In the next subsection we consider 
an example of an equivariant three-dimensional topological linear sigma
model which is a well-defined  quantum field theory by itself.

\subsection{Equivariant linear sigma model on $D\times S^1$}

In this Subsection we consider an equivariant version of the
topological linear sigma model on a non-compact three-dimensional
space $D\times S^1$ with the  target space $X=\IC^{\ell+1}$. We pick 
the flat metric \eqref{metricD} on the disk   and the
symplectic structure \eqref{symstr} on $\IC^{\ell+1}$. Let
$U_{\ell+1}$ act  on  $\IC^{\ell+1}$ via the standard
representation and $S^1$ act  on the first  factor in
$D\times S^1$ by rotations as it was introduced in Section 2. 
We would like to consider $S^1\times U_{\ell+1}$-equivariant
version of the three-dimensional topological field theory introduced
in Section 4.2. Following the reasoning of the previous Section  we
consider an  equivariance with respect to a Cartan subgroup
$G_0=S^1\times (U_1)^{\ell+1}$ of $G=S^1\times U_{\ell+1}$. The 
$G_0$-equivariant modification of the transformations
\eqref{transform} is given by \be\nonumber
\delta_{G_0}\,\varphi=\chi,\qquad
\delta_{G_0}\chi=-(\hbar\pr_{\sigma}+2\pi \imath
\beta^{-1}\pr_{\tau}+\imath \lambda)\varphi, \ee \be\nonumber
\delta_{G_0} \psi_{\zb}=F_{\zb},\qquad \delta_{G_0} F_{\zb}
=-(\hbar\pr_{\sigma}+ 2\pi \imath \beta^{-1}\pr_\tau+\imath
\lambda)\psi_{\zb}. \ee Consider the  $\delta_{G_0}$-invariant
action \eqref{Gaction3d} with $\tilde{t}=t=0$  on $N=S^1\times D$
specialized to the case of $X=\IC^{\ell+1}$ \be S_0=\int_{S^1\times
D}\, d^2z\,d\tau\,\delta_{G_0}\,\CV= \imath \int_{S^1\times D}\,
d^2z\,d\tau\,\,\Big( \pr_{\zb}\chi
\bar{\psi}_{z}+\bar{F}_{z}\pr_{\zb}\varphi +\pr_{z}\bar{\chi}
\psi_{\zb}+F_{\zb}\pr_{z}\bar{\varphi}\Big), \ee where
$$
\CV=\pr_{\zb}\varphi\,
\bar{\psi}_z+\pr_{z}\bar{\varphi}\,\psi_{\zb}.
$$
We deform the action by a $\delta_{G_0}$- and $G_0$-invariant
observable on the boundary  $\pr N=S^1\times S^1$
\be\label{Qcomthree}
S=S_0+\CO,
\ee
where
\be
\CO=\frac{\imath}{2}\beta \int_{\pr N=S^1\times
  S^1}\,d\tau\,d\sigma\,(\bar{\chi}\chi+\bar{\varphi}(
\hbar \pr_{\sigma}+2\pi \imath
\beta^{-1}\pr_{\tau}+\imath\lambda)\varphi).
 \ee
The integration over $F$ and $\psi$ localizes the functional
integral to the subspace of maps $N\to \IC^{\ell+1}$ satisfying the
equations \be\label{constr} \pr_{\zb}\varphi=0,\qquad
\pr_{z}\bar{\varphi}=0,\qquad \pr_{\zb}\chi=0,\qquad
\pr_{z}\bar{\chi}=0. \ee Let us use the expansion of the solutions
of the constraints \eqref{constr}
$$
\varphi(z,\zb,\tau)=\sum_{n=0}^{+\infty }(2\pi)^{-1/2}\varphi_n(\tau)\,z^n,
\qquad
\bar{\varphi}(z,\zb,\tau)=\sum_{n=0}^{+\infty }(2\pi)^{-1/2}\varphi_n(\tau)\,\zb^n,
$$
$$
\chi(z,\zb,\tau)=\sum_{n=0}^{+\infty }(2\pi)^{-1/2}\chi_n(\tau)\,z^n,
\qquad
\bar{\chi}(z,\zb,\tau)=\sum_{n=0}^{+\infty }(2\pi)^{-1/2}\chi_n(\tau)\,\zb^n.
$$
Then the integral factorizes into the  product of
elementary integrals \be\label{factorization}
Z=\prod_{n=0}^{\infty}Z_n, \ee where
$$Z_n=\int\,[D\varphi_n]\,[D\bar{\varphi}_n]
[D\chi_n]\,[D\bar{\chi}_n]\,\,e^{-S_n},
$$
and
$$
S_n=\frac{\imath}{2}\beta\,\int_0^{2\pi}\,d\tau\,(\bar{\chi}_n\,\chi_n+
\bar{\varphi}_n(2\pi \imath \beta^{-1}\pr_{\tau}+ \imath \hbar n+\imath \lambda)\varphi_n).
$$
Each $Z_n$ is a path integral of the form \eqref{standrrep}
and thus using previous calculations we have 
\be
Z_n=\frac{1}{1-e^{-\beta (\hbar n+ \lambda)}}.
\ee
Finally for the partition function of the three-dimensional
topological linear sigma-model on $S^1\times D$  we obtain
\be
Z=\prod_{n=0}^{+\infty}\frac{1}{1-tq^n},
\ee
where $t=e^{- \beta \lambda}$, $q=e^{- \beta \hbar}$.
This coincides with the representation given in Proposition \ref{Lemma3d}.

\section{Conclusion}

To conclude this note  we outline  some directions for future research.

In this note  we consider equivariant topological sigma models on a
disk such that the equivariance group includes  the group $S^1$ of
disk rotations. It is natural to expect that this $S^1$-equivariance
is a remnant of $Diff(S^1)$-equivariance of two-dimensional
topological quantum gravity. The relation between $S^1$-equivariance
and topological gravity is well-known. For example in  \cite{Gi1},
\cite{Gi2} the correlation functions of  $S^1$-equivariant
topological sigma models on $\IP^1$ are expressed through the
correlation functions of the topological sigma model coupled with
the topological quantum gravity. Thus one should expect that the
results of this papers can be put in the framework of a first
quantized topological string theory. However let us note that there
are also indications that  the proper interpretation of our results
should be  in terms of a second quantized topological string field theory.
The simplest soluble example of the topological string theory is
given by a pure topological gravity  completely solved
in \cite{K}. This solution can be reformulated in terms of  a
quantum field theory on a disk with a quadratic action
 playing the role of the  second quantized string theory
\cite{G}, \cite{GS}. In particular, such a formulation   
provides an intriguing analogy with the considerations
of this note and seems to deserve  further considerations.

The  construction of the functional integral representation of
local Archimedean $L$-factors uses
an integral representation of a $\Gamma$-function (see Proposition \ref{intGamma}).
Thus classical $\Gamma$-function is equal to equivariant volume of the
space of holomorphic maps  of the disk to complex plain.  
This functional integral representation should be compared with the
standard  Euler integral representation. As we
 demonstrate  in \cite{GLO5} the Euler integral representation naturally
arises as a disk partition function in the equivariant
type $B$ topological Landau-Ginzburg model on a disk
with the target space $\IC$ and the superpotential
 $W(\xi)=e^{\xi}+\lambda\xi$, $\xi\in \IC$.
This result is not surprising in view of the mirror symmetry between
the type A and type B topological sigma model (see \cite{HV} for
detailed discussion). Thus we have two integral representations of $\Gamma$-function
in terms of an infinite-dimensional equivariant symplectic volume
and  in terms of an finite-dimensional complex integral.
Taking into account the mirror symmetry
relating the two underlying topological theories,
these two integral representations  should be considered on
equal footing.

These two different integral representations of
$\Gamma$-functions are similar to  two 
constructions of the local Archimedean $L$-factors discussed in Section
1. The equivalence of
the resulting $L$-factors is a manifestation of the  local
Archimedean Langlands correspondence (see e.g. \cite{ABV}).
The analogy between mirror symmetry and local Archimedean
Langlands correspondence looks not accidental and can eventually imply that
local Archimedean Langlands correspondence follows from the  mirror
symmetry. 

In Section 4 the $q$-deformed $\Gamma$-function
was represented as a partition of an equivariant  three-dimensional
 topological linear sigma  model on
$D\times S^1$. The functional integral reduces to the functional
integral over the fields on the  boundary  $T^2=S^1\times S^1$. On the
other hand the $q$-deformed $\Gamma$-function can be identified with a
partition function of a chiral scalar field on $T^2$. This relation
between topological theory on a three-dimensional manifold and
holomorphic theory on its   boundary resembles  the relation between
conformal blocks in Wess-Zumino-Witten (WZW) theory and Chern-Simons
(CS) theory  \cite{W4}, \cite{EMSS}. Such an analog of WZW/CS
correspondence deserves considerations. 

Let us note that the proposed functional integral representation for
the classical $\Gamma$-function allows a straightforward
quantization providing  the  $q$-deformed
$\Gamma$-function. On the other hand the standard construction of
the $q$-deformed $\Gamma$-function in the classical setting is {\it
ad hoc}. One can hope that further development of our approach would
provide a canonical construction of $q$-deformations of other
classical special functions.

Local $L$-factors and their $q$-counterparts
are basic  building blocks in the description of
semi-infinite periods  associated with a type $A$ topological sigma model with the
target space $\IP^{\ell}$ and more generally
a homogeneous space of a classical Lie group
 \cite{GLO3}, \cite{GLO4}. The proposed functional integral
representations  should lead to a direct derivation of the results
of \cite{GLO4} in the framework of topological sigma models in two-
and three-dimensions.

Finally, let us stress that the main driving force of the whole
project including this note and the previous ones \cite{GLO1},
\cite{GLO2}, \cite{GLO3}, \cite{GLO4} is  to uncover the proper
geometric description of Archimedean places in arithmetic geometry.
The results of this note imply that the infinite-dimensional 
symplectic geometry could be  a proper setting  to discuss
quantum field theory models for Archimedean arithmetic geometry
seriously.

\section{Appendix:  Gaussian functional integral}

In the Appendix we describe the  standard
approach to the calculation/definition
of the Gaussian functional
integrals using $\zeta$-function
regularization \cite{RS}, \cite{H}  (see
 also  \cite{Vor}). We start with a
simple finite-dimensional Gaussian integrals
\be\label{Gausseven}
I^{\IC}_N(A)=\left(\frac{\imath}{2}\right)^{\ell+1}\,
\int_{\IC^N}\,\,e^{-\frac{1}{2}\sum_{i,j=1}^N\,\zb_iA_{ij}z_j}
\,\,\prod_{j=1}^N\,dz^j\,d\zb^j\,=
\frac{1}{\det\,A/2\pi},
\ee
where the matrix $A$ is positively defined i.e. is unitary equivalent to the
diagonal matrix with positive eigenvalues. More generally,
the Gaussian integral \eqref{Gausseven} for $A$
having complex eigenvalues $a_j$ such that ${\rm Re}(a_j)\geq 0$,
$j=1,\ldots ,N$ is  defined as a limit of the integral
 for $A$ having complex eigenvalues $a_j$ such that ${\rm Re}(a_j)>0$,
$j=1,\ldots ,N$. The resulting expression for $I^{\IC}_N(A)$
coincides with the r.h.s. of \eqref{Gausseven}.
The following integral is an example of this more
general case
\be\label{degenerate}
I=\left(\frac{\imath}{2}\right)^2\,\int_{\IC^2}\,dz d\zb\,dw\,d\wb
 \,e^{-\lambda\frac{\imath}{2}(\zb
   w+z\wb+\bar{a}w+\bar{w}a)}=
\left(\frac{2\pi}{\lambda}\right)^2,\qquad \lambda\in \IR.
\ee
Similar expression holds for Gaussian integral over anticommuting
 variables 
$$
J^{\IC}_N=\int_{\IC^{0|N}}\,e^{\imath \sum_{i,j=1}^N\,\bar{\eta}_iA_{ij}\eta_j}
\prod_{j=1}^N(\imath\,\, d\eta^j\,d\bar{\eta}^j)
=\det\,A,
$$
where we use standard  Berezin integration over anticommuting variables
$$
\int d\eta\,=0,\qquad \int \eta d\eta=1,
$$
and the sign convention in the multi-variable case is defined by 
$$
\int \eta_1\,\eta_2\cdots \eta_N\,d\eta_1 \cdots d\eta_N=1.
$$

Let  $\IC^{N|N}$ be a linear super-space
with the space of polynomial holomorphic functions \\ $\IC[z_1,\cdots
z_N,\eta_1,\cdots \eta_N]$.  We have a canonical measure on this
space
\be\label{canmes}
dm(z,\eta)=(dz_1\,d\zb_1\,dz_2\,d\zb_2 \cdots
dz_N\, d\zb_N)(d\eta_1\,d\bar{\eta}_1\,d\eta_2\,d\bar{\eta}_2 \cdots
d\eta_N\, d\bar{\eta}_N).
\ee
Then using the measure \eqref{canmes} the integral
 \eqref{Gausseven} can be rewritten as follows
\be
I^{\IC}_N(A)=\int_{\IC^{N|N}}\,dm(z,\eta)\,
e^{-\frac{1}{2}\sum_{i,j=1}^N\,\zb_iA_{ij}z_j-
\frac{\imath}{2} \sum_{j=1}^N\,\bar{\eta}_j\eta_j}=
\frac{1}{\det\,A/2\pi}.
\ee
In particular an analog of \eqref{degenerate} is given by
\be\label{degenerateodd}
\tilde{I}=\frac{1}{(2\pi)^2}\int_{\IC^{2|2}}\,dz d\zb\,dw\,d\wb\,d\eta\,d\bar{\eta}\,
d\xi\,d\bar{\xi}\,\,
 \,e^{-\lambda\frac{\imath}{2}(\zb
   w+z\wb+\bar{a}w+\bar{w}a)
+(\bar{\eta}\xi+\bar{\xi}\eta))
}=1,\qquad \lambda\in \IR.
\ee

We need certain infinite-dimensional analogs of the integrals
considered above. Instead of the matrix $A$
we have some differential operator acting in an
infinite-dimensional space of functions. To define corresponding
infinite-dimensional Gaussian integrals
one should define a notion of the determinant of the
corresponding differential operator.
Let $\CD$ be an  operator with the positive discrete
spectrum with finite multiplicities
$$
0< d_0\leq d_1\leq d_2\leq \cdots ,
$$
Consider $\zeta$-function for the operator $\CD+\lambda$
\be\label{Zetaf}
Z_{\CD}(s,\lambda)=\sum_{n=0}^{\infty}\,\frac{1}{(d_n+\lambda)^s},
\ee
where $s\in \IC$ is such that  the sum is convergent.
The sum \eqref{Zetaf}  can be continued to a meromorphic function of $s$.
The regularized determinant of the operator $\CD+\lambda$ is then
defined as follows
\be\label{zetareg}
\ln \det(\CD+\lambda)=-\pr_sZ(s,\lambda)|_{s=0}.
\ee
Consider a special case of this construction for the operator $\CD$ with
the spectrum $d_n=\rho n$, $n\in \IZ_{\geq 0}$. Corresponding
$\zeta$ -function
$$
\zeta_\rho(s,\lambda)=\sum_{n=0}^{\infty}\frac{1}{(\rho n+\lambda)^s},\qquad
  -\pi<{\rm arg}(\rho n+\lambda)\leq \pi,
$$
is basically the Hurwitz $\zeta$-functions and
has analytic continuation to all $s\in \IC/\{1\}$. We
have
$$
\zeta_\rho(0,\lambda)=\frac{1}{2}-\frac{\lambda}{\rho},\qquad
\pr_s\zeta_\rho(0,\lambda)=-\left(\frac{1}{2}-\frac{\lambda}{\rho}\right)\ln \rho+
\ln \frac{1}{\sqrt{2\pi}}\Gamma\left(\frac{\lambda}{\rho}\right).
$$
Thus for the regularized determinant of $\CD+\lambda$ we obtain
\be\label{regOne}
\Big[\prod_{n=0}^{\infty} \,(\rho n+\lambda)\Big]_{reg}=\rho ^{1/2-\lambda/\rho}
\frac{(2\pi)^{1/2}}{\Gamma(\lambda/\rho)}.
\ee
We also define the regularized determinant of the operator given by
multiplication by $\rho$ as
\be\label{regsum}
\Big[\prod_{n=0}^{+\infty}\,\rho \Big]_{reg}=\exp(\zeta(0,0)\ln \rho
)=\rho^{\frac{1}{2}}.
\ee

\begin{lem}\label{Zspec}
The following identity holds
\be\label{regTwo}
\Big[\prod_{n\in \IZ}\,(\rho n+\lambda )\Big]_{reg}\,=1-e^{2\pi
  \imath \lambda/\rho},\qquad Im(\rho)>0.
\ee
\end{lem}

\noindent {\it Proof}. We have
$$
\zeta^*_\rho(s,\lambda):=\sum_{n\in   \IZ}\,\frac{1}{(\rho n+\lambda)^s}=
\zeta_\rho(s,\lambda)+\zeta_{-\rho}(s,\lambda)-\lambda^{-s}.
$$
Simple calculation gives
$$
\frac{\pr}{\pr
  s}\zeta^*_\rho(s,\lambda)|_{s=0}=\zeta'_\rho(0,\lambda)+\zeta'_{-\rho}(0,\lambda)
+\ln \lambda
$$
$$
=-\frac{\imath \pi}{2}-\imath \pi \frac{\lambda}{\rho}+
\ln\left(\frac{1}{2\pi}\Gamma\left(\frac{\lambda}{\rho}\right)
\Gamma\left(1-\frac{\lambda}{\rho}\right)\right).
$$
Finally we obtain
$$
\Big[\prod_{n\in \IZ}\,(\rho n+\lambda)\Big]_{reg}=\exp(-\frac{\pr}{\pr
  s}\zeta^*_a(s,\lambda)|_{s=0})
=e^{\imath \pi \lambda/\rho}\frac{2\pi \imath
}{\Gamma\left(\frac{\lambda}{\rho}\right)\Gamma\left(1-\frac{\lambda}{\rho}\right)}
= (1-e^{2\pi \imath \lambda/\rho}).
$$
where we use the identity
$$
\Gamma(1-z)\Gamma(z)=\frac{\pi}{\sin (\pi z)}
$$
$\Box$

\begin{rem} Note that the regularization in  Lemma \ref{Zspec}
 is different  from a more standard regularization which uses the Fredholm definition
of the infinite determinants (see e.g. \cite{Vor}).
\end{rem}

\vskip 1cm

\noindent {\small {\bf A.G.} {\sl Institute for Theoretical and
Experimental Physics, 117259, Moscow,  Russia; \hspace{8 cm}\,
\hphantom{xxx}  \hspace{2 mm} School of Mathematics, Trinity
College, Dublin 2, Ireland; \hspace{6 cm}\hspace{5 mm}\,
\hphantom{xxx}   \hspace{2 mm} Hamilton
Mathematics Institute, Trinity College, Dublin 2, Ireland;}}

\noindent{\small {\bf D.L.} {\sl
 Institute for Theoretical and Experimental Physics,
117259, Moscow, Russia};\\
\hphantom{xxxx} {\it E-mail address}: {\tt lebedev@itep.ru}}\\

\noindent{\small {\bf S.O.} {\sl
 Institute for Theoretical and Experimental Physics,
117259, Moscow, Russia};\\
\hphantom{xxxx} {\it E-mail address}: {\tt Sergey.Oblezin@itep.ru}}

\end{document}